\documentclass[preprint,12pt]{elsarticle}
\usepackage{graphicx} % Required for inserting images
\usepackage[top=27mm, bottom=27mm, left=22mm, right=22mm]{geometry}
\usepackage{amsthm,amsmath,amssymb}
\usepackage{mathtools}
\usepackage{hyperref}
\usepackage{xcolor}
\usepackage{booktabs}
\usepackage{float}
\usepackage{bbm}
\usepackage[capitalize]{cleveref}
\newtheorem{theorem}{Theorem}[section]
\newtheorem{lemma}[theorem]{Lemma}

\newtheorem{corollary}[theorem]{Corollary}
\newtheorem{conjecture}[theorem]{Conjecture}

\theoremstyle{definition}

\newtheorem{definition}[theorem]{Definition}
\newtheorem{question}[theorem]{Question}

\newtheorem{remark}[theorem]{Remark}

\newtheorem{tthm}{Theorem}[section]

\newtheorem{llem}{Lemma}[section]
\newenvironment{llembis}[1]
  {\renewcommand{\thellem}{\ref{#1}}%
   \addtocounter{llem}{-1}%
   \begin{llem}}
  {\end{llem}}

\newcommand{\bra}[1]{\langle#1|}
\newcommand{\ket}[1]{|#1\rangle}
\newcommand{\cay}{\operatorname{Cay}}

\newcommand{\wt}{w_{\mathrm{H}}}
\newcommand{\comp}{\operatorname{comp}}
\newcommand{\per}{\operatorname{per}}

\journal{Journal of Combinatorial Theory, Series A}

\begin{document}

\begin{frontmatter}

\title{On the quantum chromatic number of Hamming and generalized Hadamard graphs}

\author[cao]{Xiwang Cao}\ead{xwcao@nuaa.edu.cn}
\author[feng]{Keqin Feng}\ead{fengkq@mail.tsinghua.edu.cn}
\author[cst-sdu]{Hexiang Huang\corref{cor1}}\ead{hexianghuang@foxmail.com}
\author[rcm-sdu]{Yulin Yang}\ead{forestyoung@mail.sdu.edu.cn}
\author[cst-sdu]{Zihao Zhang}\ead{zzhqaqhzz@gmail.com}

\affiliation[cao]{organization={School of Mathematics, Nanjing University of Aeronautics and Astronautics},             
             % addressline={},
             % city={},
             % postcode={},
             % state={},
             country={China}}
\affiliation[feng]{organization={Department of Mathematical Sciences, Tsinghua University},         
             % addressline={},
             % city={},
             % postcode={},
             % state={},
             country={China}}
\affiliation[cst-sdu]{organization={School of Cyber Science and Technology, Shandong University},             
             % addressline={},
             % city={},
             % postcode={},
             % state={},
             country={China}}
\affiliation[rcm-sdu]{organization={Research Center for Mathematics and Interdisciplinary Sciences, Shandong University},             
             % addressline={},
             % city={},
             % postcode={},
             % state={},
             country={China}}
% \affiliation[sec-ME]{organization={Key Laboratory of Cryptologic Technology and Information Security, Ministry of Education},
%              % addressline={},
%              % city={},
%              % postcode={},
%              % state={},
%              country={China}}
% \affiliation[fsc-ME]{organization={Frontiers Science Center for Nonlinear Expectations, Ministry of Education},             
%              % addressline={},
%              % city={},
%              % postcode={},
%              % state={},
%              country={China}}

\cortext[cor1]{Corresponding author.}

% \begin{abstract}
% As a fundamental metric for quantifying quantum advantage in non-local games, the quantum chromatic number reveals the power of entanglement in distributed tasks. In this paper, we investigate this parameter for $q$-ary Hamming graphs and a generalization of Hadamard graphs. Our main results establish an exponential separation between the quantum and classical chromatic numbers for both graph families, and determine the exact quantum chromatic numbers in several regimes.

% For $q$-ary Hamming graphs $H(n,q,d)$, previous bounds were largely restricted to the binary case ($q=2$) and large relative distances ($d \ge n/2$). We overcome these limitations by developing a novel linear programming framework over the Hamming scheme to construct modulus-one orthogonal representations, yielding upper bounds for arbitrary $q$ and general distances. Complementarily, we utilize the trace method to evaluate minimum eigenvalues, establishing tight spectral lower bounds for distances just below the critical threshold $(q-1)n/q$.

% For generalized Hadamard graphs over cyclic groups and finite fields, we determine their exact minimum eigenvalues. This spectral analysis reveals that the lower bound exactly matches the natural upper bound, thereby determining the precise quantum chromatic number. On the classical side, we adapt the Frankl-R\"odl forbidden intersection method to establish an exponential lower bound on the classical chromatic number, quantifying the exponential quantum advantage in these graph coloring games.
% \end{abstract}

\begin{abstract}
As a fundamental metric for quantifying quantum advantage in non-local games, the quantum chromatic number reveals the power of entanglement in distributed tasks. In this paper, we investigate this parameter for $q$-ary Hamming graphs and a generalization of Hadamard graphs. Our main results establish an exponential separation between the quantum and classical chromatic numbers for both graph families, and determine the exact quantum chromatic numbers in several regimes.
Our analysis builds on known upper and lower bounds via modulus-one orthogonal representations and minimum eigenvalues, respectively.

Previous results for Hamming graphs $H(n,q,d)$ were restricted to specific cases: the minimum eigenvalue was only identified for $d > (q-1)n/q$, while modulus-one orthogonal representations had only been constructed for the binary case ($q=2$) with $d \ge n/2$. In this work, we fill several gaps in the existing literature by developing a linear programming approach to construct modulus-one orthogonal representations for arbitrary relative distances, and using the trace method to determine the minimum eigenvalues in the regime where $d$ lies slightly below the threshold $(q-1)n/q$.

For generalized Hadamard graphs over cyclic groups and finite fields, by determining their minimum eigenvalues, we show that the spectral lower bound matches the natural upper bound on the quantum chromatic number.
On the classical side, we apply the method of forbidden intersection pattern of Frankl and R\"odl to obtain an exponential lower bound on the chromatic number, thereby quantifying the separation between the quantum and classical quantities.
\end{abstract}

% \begin{graphicalabstract}
% %\includegraphics{grabs}
% \end{graphicalabstract}

 % \begin{highlights}
 % \item There is a linear program for modulus-one orthogonal representations for Hamming graphs. 
 % \item The minimum eigenvalue of Hamming graphs with distance below and near the threshold is $K_d(2)$.
 % \item The quantum chromatic number of generalized Hadamard graphs over $\mathbb{Z}_q$ and $\mathbb{F}_q$ are determined.
 % \item Hamming graphs and generalized Hadamard graphs exhibit exponential separation between classical and quantum chromatic numbers.
 % \end{highlights}


\end{frontmatter}
\section{Introduction}

Quantum entanglement enables a fascinating phenomenon known as \textit{pseudo-telepathy}, where two spatially separated parties exhibit correlations so strong that they appear to possess telepathic powers. As described in the survey by Brassard, Broadbent, and Tapp \cite{brassard2005quantum}, this phenomenon arises in \textit{non-local games} where parties sharing an entangled quantum state can perform tasks that are classically impossible without communication.

In a general setting, two parties, Alice and Bob, receive inputs $(q_A, q_B)$ and must return answers $(a_A, a_B)$. A non-local game is called a \textit{pseudo-telepathy game} if there exists a quantum protocol that allows them to win every instance of the game with probability $1$, whereas no classical strategy (without communication) can achieve a winning probability of $1$. This strict separation highlights the power of shared entanglement. 

A canonical example of a pseudo-telepathy game is the \textit{graph coloring game} \cite{brassard2005quantum, cleve2004consequences}. Formally, this game can be cast within the general framework of \textit{quantum graph homomorphisms} \cite{MANCINSKA2016228}, which has also been specialized to characterize \textit{quantum graph isomorphisms} \cite{ATSERIAS2019289}. In the standard coloring game, a graph $G=(V, E)$ and an integer $c$ are known to both parties. A verifier sends a vertex $u$ to Alice and a vertex $v$ to Bob. Without communicating, they must return colors $c_A$ and $c_B$ from a set of $c$ colors. To win, they must satisfy the standard graph coloring constraints: they must output the same color ($c_A = c_B$) if $u = v$, and distinct colors ($c_A \neq c_B$) if $u$ and $v$ are adjacent.

Classically, Alice and Bob can win this game with probability $1$ if and only if the number of available colors $c$ is at least the classical chromatic number $\chi(G)$. If $c < \chi(G)$, any classical strategy will inevitably fail on some input pairs \cite{brassard1999cost}. 

However, utilizing the non-local properties of quantum mechanics, players can sometimes win with strictly fewer colors. The minimum number of colors required for a perfect quantum strategy is defined as the \textit{quantum chromatic number}, denoted by $\chi_Q(G)$ \cite{avis2006quantum}. 
Independently, Cameron et al. \cite{cameron2007quantum} studied this parameter (denoted therein as $\chi_q(G)$) and established its rigorous mathematical foundations.
Although the formal definition involves quantum states and measurement operators in Hilbert spaces (as detailed in Section \ref{sec:pre}, Definition \ref{def:qc}), for the purpose of outlining our results, it suffices to understand $\chi_Q(G)$ as a parameter that relaxes the conditions of classical coloring.

Much effort has been devoted to the separation between classical and quantum chromatic numbers. For small graphs, Man\v{c}inska and Roberson \cite{MancinskaRoberson2016} identified a graph on $14$ vertices exhibiting $\chi_Q(G) < \chi(G)$, which Lalonde \cite{lalonde2025quantum} later proved to be the smallest such instance.

For infinite families, the Hadamard graphs $\varOmega_n$ (where $n$ is a multiple of $4$) provide a striking example of exponential quantum advantage \cite{avis2006quantum}. In these graphs, the vertices correspond to $\pm 1$-vectors of length $n$, and two vertices are adjacent if and only if their vectors are orthogonal. They can be regarded as binary Hamming graphs $H(n, 2, n/2)$. In 1999, Brassard et al. \cite{brassard1999cost} showed that $\chi_Q(\varOmega_n) \leq n$ when $n$ is a power of $2$; this result was subsequently extended to arbitrary $n$ divisible by $4$ by Avis et al. \cite{avis2006quantum} in 2006. In contrast, Frankl and R\"odl \cite{frankl1987forbidden} demonstrated that the classical chromatic number is exponentially larger. Furthermore, when $n = 4p^k$ with $p$ being a prime, the independence number $\alpha(\varOmega_n)$, as determined by Frankl \cite{frankl1986orthogonal} and Ihringer, Tanaka \cite{ihringer2019independence}, yield an explicit lower bound for the classical chromatic number $\chi(\varOmega_n)$.

It is challenging to determine the exact value of $\chi_Q(G)$. Ji \cite{ji2013binary} proved that computing the quantum chromatic number for general graphs is NP-hard. Man\v{c}inska et al. \cite{mancinska2014graph} utilized the Lovász theta function to prove $\chi_Q(H(n,2,n/2)) = n$, a result later re-derived via spectral arguments \cite{wocjan2019spectral, mcnamara2024exact}. For a long time, Hadamard graphs were the only non-trivial family with a known quantum chromatic number. 

Very recently, Cao, Feng, and Tan \cite{cao2024quantum} determined the quantum chromatic number of another family of binary Hamming graphs, specifically $H(4t-1,2,2t)$, thereby identifying another infinite family with determined quantum chromatic numbers and exponential separation between quantum and classical chromatic number. In their work, they established an upper bound on $\chi_Q(H(n,2,d))$ for all $d \geq n/2$, while leaving the case $d < n/2$ as an open problem. 

There are rare examples of graph families exhibiting determined quantum chromatic numbers and an exponential separation from the classical chromatic numbers. This scarcity motivates the further study of quantum chromatic numbers. The following research problems are of particular interest:
\begin{itemize}
\item Construct graphs where the quantum chromatic number is strictly less than the classical chromatic number. For instance, it is known that $\chi_Q(\varOmega_n) = n$ when $4|n$, whereas $\chi(\varOmega_n) \ge (1+\varepsilon)^n$.

\item Determine the exact quantum chromatic numbers for specific classes of graphs.

\item Establish tighter lower or upper bounds for the quantum chromatic numbers of various graph classes.
\end{itemize}

% Despite the theoretical importance of Hadamard graphs, there are rare examples of graphs exhibiting an exponential separation between $\chi(G)$ and $\chi_Q(G)$. Broadening the search beyond the currently studied families is therefore crucial. In particular, deriving tight lower bounds for the classical chromatic number and upper bounds for the quantum chromatic number is {\color{red}an important method} for certifying quantum advantage. A systematic study of these bounds across a wider variety of graph classes helps to quantify the gap between classical and quantum correlations in more general settings.

To tackle these problems, we rely on two fundamental bounds on the quantum chromatic number. The first is an upper bound derived from modulus-one orthogonal representations. In this paper, graphs are assumed to be simple, namely, undirected and without multi-edges and loops. An \emph{orthogonal representation} of a graph $G$ of rank $k$ is a map $\rho: V(G) \to \mathbb{C}^k$ such that $\rho(u)$ and $\rho(v)$ are orthogonal with respect to the standard Hermitian inner product for all adjacent vertices $u,v \in V(G)$. Furthermore, a representation $\rho$ is termed \emph{modulus-one} if every entry of the vector $\rho(v)$ has unit modulus for all $v \in V(G)$. Let $\xi(G)$ and $\xi'(G)$ denote the minimum ranks of an orthogonal representation and a modulus-one orthogonal representation of $G$, respectively. Cameron et al. \cite{cameron2007quantum} established the following upper bound on the quantum chromatic number.

\begin{lemma}[{\cite[Proposition 7]{cameron2007quantum}}]
\label{lem:qc-upper}
For any graph $G$, $\chi_Q(G)\leq \xi'(G)$.
\end{lemma}

On the other hand, spectral techniques offer a lower bound. Let $\lambda_1 \ge \dots \ge \lambda_n$ be the eigenvalues of the adjacency matrix of $G$. Elphick and Wocjan \cite{elphick2019spectral} established the following Hoffman-type bound.

\begin{lemma}[\cite{elphick2019spectral}]
\label{lem:qc-lower}
Let $G$ be a graph with at least one edge. Then, 
$\chi_Q(G) \ge 1 + \frac{\lambda_1}{|\lambda_n|}$.
\end{lemma}
We remark that for regular graphs, applying this spectral lower bound reduces to determining the minimum eigenvalue of the graph.

\subsection{Our Contributions}
We investigate the quantum chromatic number and the corresponding separation properties of $q$-ary Hamming graphs, as well as a natural generalization of Hadamard graphs. Our main results establish an exponential separation between the quantum and classical chromatic numbers for both graph families, and determine the exact values in several cases.

\bigskip
\noindent\textbf{Part I: Hamming Graphs.}

The $q$-ary Hamming graph $H(n, q, d)$ is defined on the vertex set of $q$-ary $n$-tuples with adjacency determined by Hamming distance $d$. We establish tight bounds on its quantum chromatic number, making significant progress on the open question regarding the regime $d < \frac{(q-1)n}{q}$ posed in \cite{cao2024quantum}. 

\begin{theorem}[Upper Bounds of $\chi_Q(H(n,q,d))$]\label{thm:qc-upper-Hamming}
Let $n,q,d$ be positive integers with $q\ge 2$ and $d\leq n$.
\begin{enumerate}
  \item  If $d\geq\frac{(q-1)n}{q}$, then
  $$\chi_Q(H(n,q,d))\leq qd.$$
  \item  If $\frac{(q-1)n}{q} > d \geq \frac{(q-1)n-\frac{q-2}{2}-\sqrt{(q-1)n+(\frac{q-2}{2})^2}}{q}$, then
  $$\chi_Q(H(n,q,d))\leq\frac{\left(2(q-1)n-(q-2)\right)qd-q^2 d^2}{2}.$$
  \item  If $d=\delta n$ for some fixed $0<\delta<\frac{q-1}{q}$, then
  $$\chi_Q(H(n,q,d))\leq q^{H_q\left(\frac{q-1-(q-2)\delta-2\sqrt{(q-1)\delta(1-\delta)}}{q}\right)n+o(n)},$$
  where $H_q(x)=x\log_q(q-1)-x\log_qx-(1-x)\log_q(1-x)$ is the $q$-ary entropy function.
\end{enumerate}
\end{theorem}

\begin{remark}
Building on the classical results of Frankl and R\"odl \cite[Theorem 1.10]{frankl1987forbidden}, which show that $\chi(H(n,q,d))$ grows exponentially for fixed relative distances, our first two cases in \cref{thm:qc-upper-Hamming} establish an exponential separation between classical and quantum chromatic numbers.
\end{remark}

The key ingredient for establishing \cref{thm:qc-upper-Hamming} is a novel linear programming method over the \emph{Hamming scheme} \cite{delsarte1973algebraic}. This approach provides a unified framework to efficiently construct modulus-one orthogonal representations, thereby upper bounding the quantum chromatic number.
    
\begin{lemma}[Linear Programming Framework]
\label{lem:OR-LP}
The quantity $\xi'(H(n,q,d))$ is upper bounded by the value of any feasible solution to the following integer linear program:
\begin{equation*}
\begin{array}{lcl}
\text{minimize} \quad & \sum\limits_{i=0}^{n}(q-1)^i\binom{n}{i}c_i & \\[6pt]
\text{subject to} \quad\ &\begin{cases}
 \sum\limits_{i=0}^{n}c_i > 0,  \\
 \sum\limits_{i=0}^{n} c_i K_i(d) = 0, \\
 c_0,c_1,\dots,c_n \in \mathbb{Z}_{\geq 0}, 
\end{cases}
\end{array}
\end{equation*}
where $K_i(j)$ is the Krawtchouk polynomial given in \eqref{eq:krawtchouk_poly}. In fact, any feasible solution yields a modulus-one orthogonal representation.
\end{lemma}
 
To complement our upper bounds, we employ the Hoffman-type bound \cite{elphick2019spectral}, which relies on the minimum eigenvalue of the graph. For $H(n,q,d)$, the eigenvalues are given by the Krawtchouk polynomial $K_d(z)$ where $z\in\{0,1,\dots,n\}$. Van~Dam and Sotirov~\cite{van2016new} conjectured that for $d \ge \frac{(q-1)n+1}{q}$, with $d$ even when $q=2$, the minimum eigenvalue of $H(n,q,d)$ is $K_d(1)$. Alon and Sudakov~\cite{alon2000bipartite} proved this conjecture for $q=2$ when $n$ is large and $d/n$ is fixed. Dumer and Kapralova~\cite{dumer2013spherically} subsequently established the result for $q=2$ and all $n$. Finally, the conjecture was proved in full generality by Brouwer et al.~\cite{brouwer2018smallest}. While the minimum eigenvalue for $d> \frac{(q-1)n}{q}$ is $K_d(1)$, the behavior for smaller $d$ is complex. In this work, we consider the case where $d$ lies slightly below the threshold $\frac{(q-1)n}{q}$, and we identify a second threshold, above which the minimum eigenvalue of $H(n,q,d)$ is $K_d(2)$. 
    
\begin{lemma}\label{lem:eig-case2}
For $\frac{(q-1)n}{q}\geq d\geq\frac{(q-1)n-\frac{q-2}{2}-\sqrt{(q-1)n+(\frac{q-2}{2})^2-(q-1)}+1}{q}$ with $d$ even when $q=2$, and $n$ sufficiently large, the minimum eigenvalue of $H(n,q,d)$ is $K_d(2)$.
\end{lemma}

Since the minimum eigenvalue of $H(n,q,d)$ is known to be $K_d(1)$ for $d\ge \frac{(q-1)n+1}{q}$, we can combine this fact with Lemma \ref{lem:eig-case2} to obtain the following lower bounds on $\chi_Q(H(n,q,d))$.
 
\begin{theorem}[Lower Bounds of $\chi_Q(H(n,q,d))$]\label{thm:qc-lower-Hamming}
Let $n,q,d$ be positive integers with $q\ge 2$, $d\leq n$, and $d$ even when $q=2$. 
\begin{enumerate}
  \item If $d\geq\frac{(q-1)n+1}{q}$, then $$\chi_Q(H(n,q,d))\geq\dfrac{qd}{qd-(q-1)n}.$$
  \item If $\frac{(q-1)n}{q}\geq d\geq\frac{(q-1)n-\frac{q-2}{2}-\sqrt{(q-1)n+(\frac{q-2}{2})^2-(q-1)}+1}{q}$, then for sufficiently large $n$, $$\chi_Q(H(n,q,d))\geq\frac{\Big(2(q-1)n-(q-2)\Big)qd-q^2d^2}{\Big(2(q-1)n-(q-2)\Big)qd-q^2d^2-(q-1)^2(n-1)n}.$$
\end{enumerate}
\end{theorem}

\textbf{Comparison of Bounds.} Within the range of distances where the lower bounds are established (Theorem \ref{thm:qc-lower-Hamming}), the ratio between the upper and lower bounds is at most $O(n)$. This gap narrows significantly for specific distances:
\begin{itemize}
    \item At $d = \frac{(q-1)n+1}{q}$, the upper and lower bounds coincide, and thus $\chi_Q(H(n,q,d)) = (q-1)n+1$.
    \item At $d = \frac{(q-1)n}{q}$, the bounds differ by an additive constant of $q-2$. Specifically,
\[ (q-1)n-q+2 \leq \chi_Q(H(n,q,d)) \leq (q-1)n. \]
    \item For the binary case ($q=2$) at $d = \frac{n-\sqrt{n-1}+1}{2}$, the bounds differ by a multiplicative factor $\sqrt{n-1}$. Specifically, for large $n$,$$ \frac{n^2-n+2\sqrt{n-1}}{2\sqrt{n-1}} \leq \chi_Q(H(n,2,d)) \leq \frac{n^2-n+2\sqrt{n-1}}{2}. $$
\end{itemize}

% --- Part II 部分无需大改，完美衔接 ---

\bigskip
\noindent\textbf{Part II: Generalized Hadamard Graphs.}

Second, we extend our investigation to a natural generalization of the Hadamard graph over an abelian group $\mathbb{G}$ (written additively). We define the generalized Hadamard graph $\varOmega_n^{(\mathbb{G})}$ as the graph with vertex set $\mathbb{G}^n$, where two vertices $x,y\in \mathbb{G}^n$ are adjacent if and only if each element of $\mathbb{G}$ appears exactly $n / |\mathbb{G}|$ times in their difference. This construction recovers the standard Hadamard graph $\varOmega_n$ as the special case where $\mathbb{G} = \mathbb{Z}_2$.

The quantum chromatic number of generalized Hadamard graphs is always bounded above by $n$, see Lemma \ref{lem:qc-Haramard-upper}. If $\varOmega_n^{(\mathbb{G})}$ contains a clique of cardinality $n$, then we obtain a so-called
\emph{generalized Hadamard matrix} over $\mathbb{G}$, a concept first introduced by Drake
\cite{drake1979partial}. Since the seminal work of Butson \cite{butson1962generalized}, numerous constructions of generalized Hadamard matrices have been developed, see for example \cite{jungnickel1979difference}, \cite{colbourn2006difference}, and extensive surveys of diverse group structures \cite{delauney1986survey}.
When combined with bound $\omega(G)\le \chi_Q(G)$ (see \cite[Proposition 4]{cameron2007quantum}), these results determine
$\chi_Q(\varOmega_n^{(\mathbb{G})})$ for certain choices of $\mathbb{G}$ and $n$.

In this work, we focus on generalized Hadamard graphs defined over the cyclic group $\mathbb{Z}_q$ and the finite field $\mathbb{F}_q$. By employing the trace method and character sum properties, we derive the minimum eigenvalues for $\varOmega_n^{(\mathbb{Z}_q)}$ and $\varOmega_n^{(\mathbb{F}_q)}$ across a wide range of parameters. Applying the Hoffman-type bound to these results allows us to establish the lower bound $\chi_Q \ge n$ for the quantum chromatic number. Given the natural upper bound $\chi_Q \le n$, we identify new infinite families of graphs with determined quantum chromatic numbers.

% If $\varOmega_n^{(\mathbb{G})}$ contains a clique of cardinality $n$, then we obtain a so-called \emph{generalized Hadamard matrix} over $\mathbb{G}$, a concept first introduced by Drake \cite{drake1979partial}. There are many construction of generalized Hadamard matrices, combining these results with the bound $\varOmega(\varOmega_n^{(\mathbb{G})})$, we have $\varOmega(\varOmega_n^{\mathbb{G}})$ for some specific $\mathbb{G}$ and $n$. In this work, by using the trace method and properties of character sums, we respectively compute the minimum eigenvalues of $\varOmega^{(\mathbb{Z}_q)}_n$ and $\varOmega^{(\mathbb{F}_q)}_n$ for given parameters.
% Since the quantum chromatic number is naturally bounded above by $n$ (i.e., $\chi_Q \le n$), our spectral analysis implies that this upper bound is tight according to the Hoffman-type bound \cite{elphick2019spectral}. Consequently, we identify a new infinite family of graphs with explicitly known quantum chromatic numbers. 

\begin{theorem}[Exact values of quantum chromatic number]\label{thm:qc-Hadamard}
Let $q$ be a positive integer and $n$ a positive integer divisible by $q$.
\begin{enumerate}
  \item  There exists a positive integer $N = N(q)$ such that for all $n \ge N$ where $\frac{(q-1)n}{q}$ is even, we have
  \[
  \chi_Q(\varOmega_n^{(\mathbb{Z}_q)}) = n.
  \]
  \item If both $n$ and $q$ are prime powers, then
  \[
  \chi_Q(\varOmega_n^{(\mathbb{F}_q)})= n.
  \]
\end{enumerate}
\end{theorem}

In contrast to the linear growth of the quantum chromatic number, we show that the classical chromatic number grows exponentially with $n$. For the group $\mathbb{Z}_q$, we construct intersection matrices and adapt the powerful method established by Frankl and R\"odl \cite{frankl1987forbidden} to prove that $\chi(\varOmega^{(\mathbb{Z}_q)}_n)$ is exponentially large.
For the finite field $\mathbb{F}_q$, we derive a similar exponential lower bound by reducing the problem to the Hamming graph case.

\begin{theorem}[Exponentially large chromatic numbers]\label{thm:separate}
    Let $q$ be a positive integer, and let $n$ be a positive integer divisible by $q$.
\begin{enumerate}
  \item There exist positive constants $N = N(q)$ and $\varepsilon=\varepsilon(q)$ such that for all $n \ge N$ where $\frac{(q-1)n}{q}$ is even, we have 
  
  % If $\frac{(q-1)n}{q}$ is even, then there exist constants $N = N(q)\in \mathbb{N}$ and $\varepsilon=\varepsilon(q)>0$ such that for all $n \ge N$, 
  \[
  \chi(\varOmega_n^{(\mathbb{Z}_q)}) \geq (1+\varepsilon)^n.
  \]
  \item If $q$ is a prime power, there exist positive constants $N = N(q)$ and $\varepsilon=\varepsilon(q)$ such that for all $n \ge N$, with the condition that $n$ is divisible by $4$ if $q=2$, 
  \[
    \chi(\varOmega_n^{(\mathbb{F}_q)})\geq (1+\varepsilon)^n.
  \]
\end{enumerate}
\end{theorem}

Combining \cref{thm:qc-Hadamard} and \cref{thm:separate}, we confirm that these generalized Hadamard graphs also exhibit exponential quantum-classical separation.

The paper is organized as follows. In \cref{sec:pre}, we review the necessary background on quantum information and graph theory. \cref{sec:Hamming} is devoted to Hamming graphs, where we establish the upper and lower bounds given in \cref{thm:qc-upper-Hamming,thm:qc-lower-Hamming}. We then focus on generalized Hadamard graphs: \cref{sec:Hadamard} determines their exact quantum chromatic numbers (\cref{thm:qc-Hadamard}), while \cref{sec:Hadamard-sep} establishes the exponential lower bounds on their classical chromatic numbers (\cref{thm:separate}). Finally, we conclude with open problems in \cref{sec:rmk}.

\section{Preliminaries}\label{sec:pre}
For any positive integer $t$, we denote $[t] \coloneqq \{1,\dots,t\}$. \cref{subsec:quan} is intended to provide intuition regarding the quantum principles and the role of shared entanglement in the graph coloring game. Since our main results rely solely on the algebraic definition of quantum coloring (see \cref{def:qc}), this subsection may be omitted without loss of continuity.
% Throughout this paper, all graphs are assumed to be \textit{simple}, i.e., finite undirected graphs with no loops or multiple edges. 

\subsection{Basic quantum mechanics}\label{subsec:quan}
To formally define the quantum chromatic number, we first introduce the necessary notation from quantum mechanics.
We restrict our attention to finite-dimensional Hilbert spaces and use the standard Dirac bra-ket notation. 

Let $\mathbb{C}^k$ denote a column vector space and a column vector in  $\mathbb{C}^k$ is denoted by a \textit{ket} $\ket{v}$, and its conjugate transpose is denoted by a \textit{bra} $\bra{v}$. The \textit{inner product} of two vectors $\ket{v}$ and $\ket{w}$ is defined by the product $\bra{v} \ket{w}$, simply denoted by $\langle v|w\rangle$.
 In quantum information, the state of a system (or simply a particle) is treated as a unit vector in $\mathbb{C}^k$ (i.e., a vector $\ket{v}$ satisfying $\langle v|v\rangle=1$). Such a system is usually referred to as a \textit{qubit} when $k=2$.

\begin{definition}\label{def:complete orthogonal system}
A set of matrices $\{P_1, P_2, \ldots, P_t\}$ acting on $\mathbb{C}^k$ is called a \textit{complete orthogonal system} if:
\begin{enumerate}
    \item each $P_i$ is idempotent, i.e., $P_i^2 = P_i$;
    \item the system is mutually orthogonal, i.e., $P_i P_j = {0}$ for any $i \neq j$;
    \item the system is complete, i.e., $\sum_{i \in [t]} P_i = I_k$.
\end{enumerate}
\end{definition}
Suppose $\mathbb{C}^k$ admits a direct sum decomposition $\mathbb{C}^k = V_1 \oplus V_2 \oplus \cdots \oplus V_t$. Let $P_i$ be the projective operator from $\mathbb{C}^k$ onto $V_i$ along this decomposition, defined by $P_i(v_1 + \cdots + v_t) = v_i$ for all $v_j \in V_j$ and $j \in [t]$. Then, the set $\{P_i : i \in [t]\}$ forms a complete orthogonal system. Conversely, given a complete orthogonal system acting on $\mathbb{C}^k$, we obtain a direct sum decomposition $\mathbb{C}^k = \mathrm{im}(P_1) \oplus \cdots \oplus \mathrm{im}(P_t)$. Specifically, an \textit{orthogonal} decomposition of $\mathbb{C}^k$ corresponds to a complete orthogonal system of \textit{Hermitian} projectors (in the context of quantum information, a Hermitian idempotent is usually called an orthogonal projector).

Assume an orthogonal decomposition $\mathbb{C}^k = V_1 \perp \cdots \perp V_t$, corresponding to a complete set of orthogonal projectors $\{P_i : i \in [t]\}$. Consider a state $\ket{v} \in \mathbb{C}^k$ (a unit vector). The state $\ket{v}$ can be uniquely expressed as $\ket{v} = \ket{v_1} + \cdots + \ket{v_t}$ with $\ket{v_i} \in V_i$. When a measurement defined by this decomposition (or equivalently, by the system $\{P_i\}$) is performed, the state $\ket{v}$ collapses to the normalized vector $\frac{\ket{v_i}}{\|\ket{v_i}\|}$ with probability $\|\ket{v_i}\|^2$. Observing the outcome $i$ corresponds to observing that the state has collapsed into the subspace $V_i$. The probability distribution is valid since $\sum_{i \in [t]} \|\ket{v_i}\|^2 = \|\ket{v}\|^2 = 1$. This process is known as a \textit{projective measurement}.

Using the relation $\ket{v_i} = P_i\ket{v}$, we can express the probability and the post-measurement state in terms of the projectors $P_i$:
\begin{align*}
    \Pr(i) &= \|\ket{v_i}\|^2 = \langle v_i|v_i\rangle = \bra{v}P_i^\dagger P_i\ket{v} = \bra{v}P_i^2\ket{v} = \bra{v}P_i\ket{v}, \\
    \text{Post-state} &= \frac{\ket{v_i}}{\|\ket{v_i}\|} = \frac{P_i\ket{v}}{\sqrt{\bra{v}P_i\ket{v}}}.
\end{align*}
For simplicity, we refer to the set $\{P_i : i \in [t]\}$ as a projective measurement. As background, projective measurements are sufficient to understand the quantum graph coloring, so we focus on them. (While the general framework of quantum operations encompasses POVMs and unitary evolutions, these lie beyond the scope of our current discussion.)

Suppose $k$ is a positive integer. For each $i \in [k]$, let $\ket{i} \in \mathbb{C}^k$ denote the column vector with a $1$ in the $i$-th position and $0$ elsewhere. The orthonormal basis $\{\ket{i} : i \in [k]\}$ is referred to as the \textit{computational basis}. The corresponding set of projectors $\{\ket{i}\bra{i} : i \in [k]\}$ serves as a fundamental example of a projective measurement. Note that in a general projective measurement, the subspaces $V_i$ need not have the same dimension.

Consider performing two sequential measurements on a quantum state $\ket{v} \in \mathbb{C}^k$, defined first by the set $\{P_i : i \in [t]\}$ and subsequently by $\{Q_j : j \in [s]\}$. This procedure yields an outcome pair $(i,j) \in [t] \times [s]$. To determine the joint probability $\Pr(i,j)$, we analyze the process step by step. The probability of obtaining the first outcome $i$ is $\bra{v}P_i\ket{v}$. Upon this outcome, the state $\ket{v}$ collapses to the normalized vector $\ket{v'} = \frac{P_i\ket{v}}{\sqrt{\bra{v}P_i\ket{v}}}$. According to the postulates of quantum mechanics, the conditional probability of obtaining the second outcome $j$ given $i$ is:
\[
    \Pr(j|i) = \bra{v'}Q_j\ket{v'} = \frac{1}{\bra{v}P_i\ket{v}} \bra{v}P_i^\dagger Q_j P_i\ket{v} = \frac{\bra{v}P_i Q_j P_i\ket{v}}{\bra{v}P_i\ket{v}}.
\]
Consequently, the joint probability is given by $\Pr(i,j) = \bra{v}P_i Q_j P_i\ket{v}$. This formula holds generally: if $\bra{v}P_i\ket{v} = 0$, both the probability $\Pr(i)$ and the joint probability vanish (consistent with $P_i\ket{v} = {0}$); otherwise, $\Pr(i,j) = \Pr(i)\Pr(j|i)$, which also yields the expression above.

If the two measurements commute (i.e., $P_i Q_j = Q_j P_i$ for all $i,j$), the expression simplifies to:
\[
    \Pr(i,j) = \bra{v}P_i Q_j \ket{v}.
\]
In this commuting case, the order in which the measurements $\{P_i\}$ and $\{Q_j\}$ are applied does not affect the joint probability distribution.

Two quantum systems can be combined to form a composite system. For instance, the tensor product of the spaces $\mathbb{C}^k$ and $\mathbb{C}^\ell$ yields the joint space $\mathbb{C}^k \otimes \mathbb{C}^\ell$. A bipartite state $\ket{v} \in \mathbb{C}^k \otimes \mathbb{C}^\ell$ describes a system of two particles that can be distributed to two spatially separated observers, Alice and Bob. Mathematically, this separation implies that Alice's measurements are of the form $\{P_i \otimes I_{\ell} : i \in [t]\}$, while Bob's are of the form $\{I_k \otimes Q_j : j \in [s]\}$. Since operators acting on disjoint subsystems commute, the sequential measurement principle derived above applies, implying that the order of measurement is irrelevant. The resulting joint probability is given by
\[\Pr(i,j) = \bra{v} P_i \otimes Q_j \ket{v}.\]
For convenience, we denote Alice's and Bob's measurements as $\{P_i\}$ and $\{Q_j \}$, respectively.

\subsection{Quantum coloring for graphs}
We first formalize the notion of quantum coloring. As rigorously formulated by Cameron et al.~\cite{cameron2007quantum}, the quantum chromatic number $\chi_Q$ relates to the existence of a specific system of projectors.

\begin{definition}\label{def:qc}
    Let $G$ be a simple graph and $r$ be a positive integer. Consider a set of Hermitian matrices
    \[
        \left\{P_{v,\alpha} : v \in V(G), \alpha \in [r]\right\} \subseteq \mathbb{C}^{k \times k},
    \]
    where $[r] = \{1,2,\ldots,r\}$ and $k$ is a positive integer. The set $\{P_{v,\alpha} : v \in V(G), \alpha \in [r]\}$ is called a \emph{quantum $r$-coloring} of $G$ if it satisfies the following conditions:
    \begin{itemize}
        \item[(1)] For every vertex $v \in V(G)$, the set $\{P_{v,\alpha} : \alpha \in [r]\}$ forms a complete orthogonal system (see \cref{def:complete orthogonal system}).
        \item[(2)] For any color $\alpha \in [r]$ and any two adjacent vertices $u, v \in V(G)$, we have the orthogonality condition $P_{u,\alpha} P_{v,\alpha} = O$.
    \end{itemize}
\end{definition}

A quantum $r$-coloring $\{P_{v,\alpha}\}_{v,\alpha}$ provides a perfect strategy for the graph coloring game using shared entanglement. Specifically, Alice and Bob share a maximally entangled state $\ket{v} \in \mathbb{C}^{k} \otimes \mathbb{C}^k$, defined as
\[
    \ket{v} = \frac{1}{\sqrt{k}}\sum_{i=1}^k \ket{i} \otimes \ket{i}.
\]
The two particles of $\ket{v}$ are distributed to Alice and Bob, respectively. Upon receiving vertices $x$ and $y$ respectively, Alice performs the measurement $\{P_{x,i}\}_{i}$ and Bob performs the measurement $\{\overline{P_{y,j}}\}_{j}$ (where $\overline{P}$ denotes the complex conjugate) on their own particles. Alice then reports her measurement result $\alpha$ to the verifier, and Bob similarly reports his measurement result $\beta$.

The probability that Alice reports color $\alpha$ and Bob reports color $\beta$ is given by
\begin{align*}
    \Pr(\alpha,\beta) &= \bra{v} (P_{x,\alpha} \otimes \overline{P_{y,\beta}}) \ket{v} \\
    &= \frac{1}{k} \left( \sum_{i \in [k]} \bra{i} \otimes \bra{i} \right) (P_{x,\alpha} \otimes \overline{P_{y,\beta}}) \left( \sum_{j \in [k]} \ket{j} \otimes \ket{j} \right) \\
    &= \frac{1}{k} \sum_{i,j \in [k]} (\bra{i} P_{x,\alpha} \ket{j}) (\bra{i} \overline{P_{y,\beta}} \ket{j}) \\
    &= \frac{1}{k} \sum_{i,j \in [k]} \bra{i} P_{x,\alpha} \ket{j} \bra{j} P_{y,\beta} \ket{i} \quad \text{(since $P_{y,\beta}$ is Hermitian)} \\
    &= \frac{1}{k} \sum_{i \in [k]} \bra{i} P_{x,\alpha} \left( \sum_{j \in [k]} \ket{j} \bra{j} \right) P_{y,\beta} \ket{i} \\
    &= \frac{1}{k} \sum_{i \in [k]} \bra{i} P_{x,\alpha} P_{y,\beta} \ket{i} = \frac{1}{k} \mathrm{Tr}(P_{x,\alpha} P_{y,\beta}).
\end{align*}

The conditions in \cref{def:qc} ensure that this strategy wins with probability $1$. First, if the vertices $x$ and $y$ are adjacent ($x \sim y$), the probability of obtaining the same color $\alpha = \beta$ is
\[
    \Pr(\alpha,\alpha) = \frac{1}{k}\mathrm{Tr}(P_{x,\alpha}P_{y,\alpha}) = 0,
\]
which follows from the orthogonality condition $P_{x,\alpha}P_{y,\alpha} = O$.
Second, if the vertices are the same ($x=y$), the probability of obtaining different colors $\alpha \neq \beta$ is
\[
    \Pr(\alpha,\beta) = \frac{1}{k}\mathrm{Tr}(P_{x,\alpha}P_{x,\beta}) = 0,
\]
since $\{P_{x,\gamma}\}_\gamma$ forms an orthogonal system. 

Cameron et al. \cite{cameron2007quantum} showed that a quantum protocol employing general measurements and shared entanglement used to win the graph coloring game turns out to be a strategy utilizing projective measurements and maximal entangled state as described above. Building on this equivalence, the quantum chromatic number of a graph $G$, denoted by $\chi_Q(G)$, can be defined as the minimum integer $r$ for which a quantum $r$-coloring of $G$ exists.

For the reader's convenience, we restate the following lemmas, which link the quantum chromatic number to the orthogonal rank and the spectrum of graphs, respectively.

\begin{llembis}{lem:qc-upper}[{\cite[Proposition 7]{cameron2007quantum}}]
    For any graph $G$, $\chi_Q(G)\leq \xi'(G)$.
\end{llembis}

\begin{proof}
Let $k = \xi'(G)$. Suppose $\rho: V(G)\rightarrow \mathbb{C}^k$ is a modulus-one orthogonal representation of $G$.

Denote $\zeta_k = e^{\frac{2\pi\sqrt{-1} }{k}}$. For each $x\in V(G)$ and $\alpha\in[k]$, define
$P_{x,\alpha}=\Lambda_x\ket{\alpha}\bra{\alpha}\Lambda_x^{-1}$, 
where $\Lambda_x = \operatorname{diag}(\rho(x)_i:i\in [k])$ and
$$\ket{\alpha}=\frac{1}{\sqrt{k}}
\begin{bmatrix}
\zeta_k^{\alpha0}\\ \zeta_k^{\alpha1}\\ \vdots \\ \zeta_k^{\alpha(k-1)}\end{bmatrix}.$$
Clearly, each $P_{x,\alpha}$ is Hermitian and $\{P_{x,\alpha}:\alpha\in [k]\}$ forms a complete orthogonal system for each $x\in V(G)$. It remains to verify the orthogonality condition, namely, $P_{x,\alpha}P_{y,\alpha} = 0$ for all adjacent $x,y\in V(G)$, and $\alpha\in [k]$. We have
\[P_{x,\alpha}P_{y,\alpha} = (\Lambda_x\ket{\alpha}\bra{\alpha}\Lambda_x^{-1})(\Lambda_y\ket{\alpha}\bra{\alpha}\Lambda_y^{-1})=0,\]
since 
\[\bra{\alpha}\Lambda_x^\dagger\Lambda_y\ket{\alpha} = 
\begin{bmatrix}
    \dots\ \zeta_k^{-\alpha i}\ \dots
\end{bmatrix}
\begin{bmatrix}
   \ddots & & \\
    & \rho(x)_i^* & \\
    & & \ddots
\end{bmatrix}
\begin{bmatrix}
   \ddots & & \\
    & \rho(y)_i & \\
    & & \ddots
\end{bmatrix}
\begin{bmatrix}
    \vdots\\
    \zeta_k^{\alpha i}\\
    \vdots
\end{bmatrix}
= \sum_{i = 0}^{k-1} \rho(x)_i^*\rho(y)_i = 0.
\]
Therefore, $\{P_{x,\alpha}:x\in V(G), \alpha\in [k]\}$ is a quantum $k$-coloring of $G$, hence $\chi_Q(G)\le k= \xi'(G)$.
\end{proof}

\begin{llembis}{lem:qc-lower}[\cite{elphick2019spectral}]
Let $G$ be a graph with at least one edge. Then, $\chi_Q(G) \ge 1 +\frac{\lambda_1}{|\lambda_n|}$.
\end{llembis}

% We next recall two fundamental bounds for the quantum chromatic number: the first is an upper bound based on modulus-one orthogonal representations.

% An \emph{orthogonal representation} of a graph $G$ with dimension $K$ is a map 
% $\rho: V(G) \to \mathbb{C}^K$ such that $\rho(u)$ and $\rho(v)$ are orthogonal 
% with respect to the complex inner product for all adjacent vertices $u,v \in V(G)$. 
% Moreover, the representation $\rho$ is called \emph{modulus-one} if all coordinates of $\rho(v)$ 
% have modulus one for every $v \in V(G)$. Let $\xi(G)$ and $\xi'(G)$ denote the minimum dimension of an orthogonal representation and of a modulus-one orthogonal representation, respectively. 

% As established by Cameron et al. in \cite{cameron2007quantum}, the quantum chromatic number satisfies the following upper bound:
% \begin{lemma}\label{lem:qc-upper}
% For any graph $G$, $\chi_Q(G)\leq \xi'(G)$.
% \end{lemma}

% Complementarily, spectral techniques provide a lower bound for $\chi_Q(G)$. 
% Let $\lambda_1 \ge \lambda_2 \ge \dots \ge \lambda_n$ denote the eigenvalues of the adjacency matrix of $G$.
% The following Hoffman-type bound was established in \cite{elphick2019spectral}:
% \begin{lemma}\label{lem:qc-lower}
% Let $G$ be a graph with at least one edge. Then, 
% $\chi_Q(G) \ge 1 - \frac{\lambda_1}{\lambda_n}$.
% \end{lemma}
% In particular, if $G$ is $d$-regular, then $\lambda_1 = d$, and the bound becomes $1 - d/\lambda_n$.

\subsection{Krawtchouk polynomials}

We first recall some standard graph parameters. Given two graphs $G$ and $H$, a \emph{homomorphism} from $G$ to $H$ is a map $\phi : V(G) \to V(H)$ such that if $u$ and $v$ are adjacent in $G$, then $\phi(u)$ and $\phi(v)$ are adjacent in $H$. The chromatic number of $G$, denoted by $\chi(G)$, is the minimum positive integer $r$ such that there exists a homomorphism from $G$ to the complete graph on $r$ vertices. The clique number of $G$, denoted by $\varOmega(G)$, is the maximum positive integer $r$ such that there exists a homomorphism from the complete graph on $r$ vertices to $G$. The independence number of $G$, denoted by $\alpha(G)$, is the maximum size of an independent set in $G$ and is equal to the clique number of the complement of $G$. It is noteworthy that $\chi(G)\alpha(G) \geq |V(G)|$.

Let $\mathbb{G}$ be a group and $S \subseteq \mathbb{G}$ a subset. The \emph{Cayley graph} $\cay(\mathbb{G}, S)$ is the directed graph with vertex set $\mathbb{G}$, where a directed edge $(g,h)$ exists if and only if $g^{-1}h \in S$. If $S$ is inverse-closed (i.e., $S^{-1} = S$), then $\cay(\mathbb{G},S)$ is undirected.

Both Hamming and composition graphs can be realized as Cayley graphs over abelian groups. 
Let $\varSigma$ be a $q$-ary alphabet endowed with an additive group structure $\mathbb{G}=(\varSigma,+)$.
For a tuple $x \in \mathbb{G}^n$, the \emph{Hamming weight} $\wt(x)$ is the number of coordinates not equal to $0_\mathbb{G}$, and the \emph{composition} $\comp(x)$ is a tuple $\boldsymbol{i}=(i_g : g \in \mathbb{G})$ where $i_g$ is the number of coordinates of $x$ equal to $g$. 
The Hamming graph $H(n,q,d)$ and the composition graph $H_C(n,\mathbb{G},\boldsymbol{i})$ are isomorphic to the Cayley graphs $\cay(\mathbb{G}^n, S_d)$ and $\cay(\mathbb{G}^n, S_{\boldsymbol{i}})$, respectively, where 
\[
S_d = \{x \in \mathbb{G}^n : \wt(x) = d\} \text{~~and~~} S_{\boldsymbol{i}} = \{x \in \mathbb{G}^n : \comp(x) = \boldsymbol{i}\}.
\] 
In particular, when $q \mid n$, we denote the balanced composition $(n/q, \dots, n/q)$ by $\boldsymbol{n/q}$. Note that the generalized Hadamard graph over $\mathbb{G}$ (with $|\mathbb{G}|=q$) of length $n$ is defined as $\varOmega_n^{(\mathbb{G})}=H_C(n,\mathbb{G},\boldsymbol{n/q})$, which is therefore isomorphic to $\cay(\mathbb{G}^n, S_{\boldsymbol{n/q}})$.

Recall that both Hamming and Composition graphs are Cayley graphs over an abelian group. 
The eigenvalues of an abelian Cayley graph can be expressed elegantly. Before presenting their eigenvalues, we first recall some basics of the characters of an abelian group.

Let $\mathbb{G}$ be a finite abelian group. A character $\varphi$ of $\mathbb{G}$ is a homomorphism 
from $\mathbb{G}$ to the multiplicative group of complex numbers, i.e., $\varphi:\mathbb{G}\to\mathbb{C}^{\times}$, 
satisfying $\varphi(xy)=\varphi(x)\varphi(y)$. 
Let $\widehat{\mathbb{G}}$ denote the set of all characters of $\mathbb{G}$. 
For each $x \in \mathbb{G}$, we have $\varphi(x)^{|\mathbb{G}|} = \varphi(|\mathbb{G}|x) = \varphi(0_\mathbb{G}) = 1$. 
Thus the values of $\varphi$ are $|\mathbb{G}|$-th roots of unity. 
Moreover, $\widehat{\mathbb{G}}$ forms an abelian group under entrywise multiplication, and there is a canonical isomorphism from $\mathbb{G}$ to $\widehat{\mathbb{G}}$. 
Under this isomorphism, we can relabel $\widehat{\mathbb{G}} = \{\varphi_g : g \in \mathbb{G}\}$ such that $\varphi_{g+h} = \varphi_g \varphi_h$ and $\varphi_g(h) = \varphi_h(g)$ for all $g,h \in \mathbb{G}$.

It is well known that for $\mathbb{G} = \mathbb{Z}_q^n$, 
the character group is given by $\widehat{\mathbb{Z}_q^n} = \{\varphi_x : x \in \mathbb{Z}_q^n\}$, 
with $\varphi_x(y) = \zeta_q^{x \cdot y}$ for all $y \in \mathbb{Z}_q^n$.
Similarly, for $\mathbb{G} = \mathbb{F}_q^n$, where $q$ is a power of a prime $p$, 
we have $\widehat{\mathbb{F}_q^n} = \{\psi_\alpha : \alpha \in \mathbb{F}_q^n\}$, 
and the characters are defined by $\psi_\alpha(\beta) = \zeta_p^{\operatorname{Tr}(\alpha \cdot \beta)}$ for all $\beta \in \mathbb{F}_q^n$, where $\operatorname{Tr}(\cdot)$ denote the trace function from $\mathbb{F}_q$ to $\mathbb{F}_p$.
In both cases, the dot product $x \cdot y$ is defined as $\sum_{i=1}^n x_i y_i$, and $\zeta_r$ is a primitive $r$-th root of unity.

We now recall a classical result on the spectrum of abelian Cayley graphs using characters.
\begin{lemma}[\cite{lovasz1975spectra}]
\label{thm:eigcay}
    Let $\mathbb{G}$ be a finite abelian group, and let $S$ be a subset of $\mathbb{G}$. Then the eigenvalues of the Cayley graph $\cay(\mathbb{G},S)$ are given by
    \begin{equation*}
        \lambda_g=\sum_{s\in S}\varphi_g(s), \quad g\in{\mathbb{G}}.
    \end{equation*}
    Moreover, the eigenvector corresponding to $\lambda_\varphi$ is $(\varphi(x):x\in \mathbb{G})^\top$.
\end{lemma}

As a consequence, this result provides an explicit expression for the eigenvalues of Hamming and composition graphs.

For Hamming graphs, the eigenvalues of $H(n,q,i)\cong\cay(\mathbb{G}^n,S_i)$ are given by
\begin{equation*}
    \lambda_x = \sum_{y \in S_i} \varphi_x(y), \quad x \in \mathbb{G}^n.
\end{equation*}
A direct calculation in \cite{delsarte1973algebraic} shows that $\lambda_x$ depends only on the Hamming weight of $x$. 
More precisely, if $\wt(x) = j$, then $\lambda_x$ equals the degree-$i$ Krawtchouk polynomial evaluated at $j$, denoted by $K_i(j)$, which is defined by
\begin{equation}
\label{eq:krawtchouk_poly}
    K_i(j) = \sum_{k=0}^i (-1)^k (q-1)^{i-k} \binom{j}{k} \binom{n-j}{i-k}.
\end{equation}
This fact leads naturally to the following reciprocal property, which a simple double-counting argument can derive:
\begin{equation*}
    (q-1)^j\binom{n}{j}K_{i}(j)=\sum_{x\in S_j}\sum_{y\in S_i}\varphi_x(y)=\sum_{y\in S_i}\sum_{x\in S_j}\varphi_y(x)= (q-1)^i\binom{n}{i}K_{j}(i).
\end{equation*}

For composition graphs, the eigenvalues of $H(n,\mathbb{G},\boldsymbol{i}) \cong \cay(\mathbb{G}^n,S_{\boldsymbol{i}})$ are given by
\begin{equation*}
\mu_x = \sum_{y \in S_{\boldsymbol{i}}} \varphi_x(y), \quad x \in \mathbb{G}^n.
\end{equation*}
Let $\boldsymbol{z} = (z_g : g \in \mathbb{G})$ be a tuple of indeterminates, and denote the monomial 
$\prod_{g \in \mathbb{G}} z_g^{i_g}$ by $\boldsymbol{z}^{\boldsymbol{i}}$. 
We use the bracket notation $[\cdot]$ to denote the coefficient extraction operator; that is, 
$[\boldsymbol{z}^{\boldsymbol{i}}] f(\boldsymbol{z})$ gives the coefficient of $\boldsymbol{z}^{\boldsymbol{i}}$ in the expansion of $f(\boldsymbol{z})$.

To encode the above character sums over $n$-tuples in $S_{\boldsymbol{i}}$, we consider the generating function
\[
\prod_{k=1}^n \left( \sum_{g \in \mathbb{G}} \varphi_{x_k}(g) z_g \right),
\]
where for each $k$, $\varphi_{x_k}$ is the character of $\mathbb{G}$ corresponding to the $k$-th coordinate of $x$.  

Expanding the product, each monomial corresponds to an $n$-tuple $y \in \mathbb{G}^n$, with the exponent of $z_g$ equal to the number of coordinates of $y$ that are equal to $g$.  
It follows that the coefficient of $\boldsymbol{z}^{\boldsymbol{i}}$ in this expansion is precisely the sum of the product of characters over all $y \in S_{\boldsymbol{i}}$, that is,
\[
\sum_{y \in S_{\boldsymbol{i}}} \varphi_x(y)=\sum_{y \in S_{\boldsymbol{i}}} \prod_{k=1}^n \varphi_{x_k}(y_k)
= [\boldsymbol{z}^{\boldsymbol{i}}] \prod_{k=1}^n \left( \sum_{g \in \mathbb{G}} \varphi_{x_k}(g) z_g \right).
\]

Another expression for $\mu_x$ shows that its value is determined solely by the composition of $x$. Specifically,  if $\mathrm{comp}(x) = \boldsymbol{j} = (j_g : g \in \mathbb{G})$, then
\begin{equation*}
\mu_x
= [\boldsymbol{z}^{\boldsymbol{i}}]
\prod_{k=1}^n\left(\sum_{g\in \mathbb{G}}\varphi_{x_k}(g)z_g\right)
= [\boldsymbol{z}^{\boldsymbol{i}}]
\prod_{h\in \mathbb{G}}\left(\sum_{g\in \mathbb{G}}\varphi_{h}(g)z_g\right)^{j_h}.
\end{equation*}

We define the \emph{generalized Krawtchouk polynomial} over $\mathbb{G}$ \footnote{In the special case $\mathbb{G} = \mathbb{F}_q$, this polynomial coincides with the one studied by Sookoo \cite{sookoo2020generalized}.}  as
\begin{equation*}
K^{(\mathbb{G})}_{\boldsymbol{i}}(\boldsymbol{j})
= [\boldsymbol{z}^{\boldsymbol{i}}]
\prod_{h\in \mathbb{G}}\left(\sum_{g\in \mathbb{G}} \varphi_{h}(g) z_g \right)^{j_h},
\end{equation*}

Finally, by a simple double-counting argument, the generalized Krawtchouk polynomials also satisfy the reciprocal law:
\begin{equation*}
\binom{n}{\boldsymbol{j}} K^{(\mathbb{G})}_{\boldsymbol{i}}(\boldsymbol{j})
= \sum_{x \in S_{\boldsymbol{j}}} \sum_{y \in S_{\boldsymbol{i}}} \varphi_x(y)
= \sum_{y \in S_{\boldsymbol{i}}} \sum_{x \in S_{\boldsymbol{j}}} \varphi_y(x)
= \binom{n}{\boldsymbol{i}} K^{(\mathbb{G})}_{\boldsymbol{j}}(\boldsymbol{i}).
\end{equation*}
where we use the notation 
$\binom{n}{\boldsymbol{i}}$ for the multinomial coefficient 
$\frac{n!}{\prod_{g \in \mathbb{G}} i_g!}$.

\section{Quantum chromatic number of some Hamming graphs}\label{sec:Hamming}
 The goal of this section is to present the proofs of \cref{thm:qc-upper-Hamming,thm:qc-lower-Hamming}.
\subsection{A linear programming approach to orthogonal representations}\label{sec3:LP bound}
In this subsection, we derive upper bounds on $\xi'(H(n,q,d))$ using a linear programming approach. We begin by recalling some basic facts about the Hamming scheme. 

Let $A_i$ be the adjacency matrix of $H(n,q,i)$ for $i = 0,1,\dots,n$.
Then $A_0 = I, A_1, \dots, A_n$ span the \emph{Bose–Mesner algebra} of the Hamming scheme, denoted by
\begin{equation*}
\mathbb{A} = \operatorname{span}_{\mathbb{C}}\{A_0, A_1, \ldots, A_n\}.
\end{equation*}

Denote $\ket{\varphi_x} = (\varphi_x(y) : y \in \mathbb{G}^n)^\top$.
As shown in the previous section, we have
\begin{equation*}
A_i \cdot \ket{\varphi_x} = K_i(\wt(x)) \cdot \ket{\varphi_x}, \quad i = 0,1,\dots,n.
\end{equation*}
This implies that all matrices in $\mathbb{A}$ can be simultaneously diagonalized.
Consequently, viewing $\mathbb{C}^{q^n}$ as an $\mathbb{A}$-module, it decomposes orthogonally into a direct sum of $\mathbb{A}$-submodules:
\begin{equation*}
\mathbb{C}^{q^n} = V_0 \perp V_1 \perp \dots \perp V_n,
\end{equation*}
where $V_j = \operatorname{span}_{\mathbb{C}}\{\ket{\varphi_x} : \wt(x) = j\}$.
Let $E_j$ be the orthogonal projection from $\mathbb{C}^{q^n}$ onto $V_j$, namely
\begin{equation*}
    E_j = \sum_{x\in S_j}\ket{\varphi_x}\bra{\varphi_x}, \quad j=0,1,\dots, n, 
\end{equation*}
where $\bra{\varphi_x}$ denote the conjugate transpose of $\ket{\varphi_x}$.
Then $E_0=J, E_1, \dots, E_n$ form a complete orthogonal system, and each adjacency matrix $A_i$ can be expressed as a linear combination of these projectors, with the corresponding eigenvalues as coefficients:
\begin{equation*}
\label{eq:linear_relation_between_Ai_and_Ej}
    A_i = K_i(0)E_0+K_i(1)E_1+\dots+K_i(n)E_n, \quad i = 0,1,\dots,n.
\end{equation*}
Therefore, $\mathbb{A}$ lies in the span of $E_0, E_1, \dots, E_n$.  
Since $A_0, A_1, \dots, A_n$ are clearly linearly independent, we have $\dim(\mathbb{A}) = n+1$.  
By comparing dimensions, it follows that the span of $E_0, E_1, \dots, E_n$ is exactly $\mathbb{A}$.  
It is well known (see, e.g., \cite{delsarte1973algebraic}) that the following change-of-basis formula holds:
\begin{equation*}
        E_j = \frac{1}{q^n}\bigl(K_j(0)A_0 + K_j(1)A_1 + \cdots + K_j(n)A_n\bigr), \quad j=0,1,\dots,n.  
\end{equation*}

The following lemma presents an explicit construction of a modulus-one orthogonal representation of Hamming graphs. For the reader's convenience,  we restate it here from the introduction.

\begin{llembis}{lem:OR-LP}
    The quantity $\xi'(H(n,q,d))$ is upper bounded by the objective value of any feasible solution to the following integer linear program:
    \begin{equation*}
    \begin{array}{lcl}
    \text{minimize} \quad & \sum\limits_{i=0}^{n}(q-1)^i\binom{n}{i}c_i & \\[6pt]
    \text{subject to} \quad\ &\begin{cases}
     \sum\limits_{i=0}^{n}c_i > 0,  \\
     \sum\limits_{i=0}^{n} c_i K_i(d) = 0, \\
     c_0,c_1,\dots,c_n \in \mathbb{Z}_{\geq 0}, 
    \end{cases}
    \end{array}
    \end{equation*}
    where $K_i(j) = [z^i](1-z)^j(1+(q-1)z)^{n-j}$ is the $q$-ary Krawtchouk polynomial. In fact, any feasible solution yields a modulus-one orthogonal representation.
\end{llembis}

% \begin{lemma}
% \label{lem:OR-LP}
% The quantity $\xi'(H(n,q,d))$ is upper bounded by the value of any feasible solution to the following linear program:
% \begin{equation*}
% \begin{array}{lcl}
% \text{minimize} \quad & \sum\limits_{i=0}^{n}(q-1)^i\binom{n}{i}c_i & \\[6pt]
% \text{subject to} \quad 
% &  \sum\limits_{i=0}^{n}c_i > 0, & \\
% &  \sum\limits_{i=0}^{n} c_i K_i(d) = 0, & \\
% & c_0,c_1,\dots,c_n \in \mathbb{Z}_{\geq 0}. &
% \end{array}
% \end{equation*}
% % Moreover, any feasible solution $(x_0,\dots,x_n)$ yields an orthogonal representation of $H(n,q,d)$ in dimension $\sum_{i=0}^{n}(q-1)^i\binom{n}{j}x_i$, mapping each vertex $a\in\mathbb{G}^n$ toto the $a$-th row of the matrix obtained by concatenating $x_j$ copies of $R_i$ for all $i$ with $x_i\ne 0$, where $R_i$ denote the matrix whose columns are the vectors $\ket{\varphi_a}$ for $a \in S_i$.
% \end{lemma}

\begin{proof}
Let $(c_0,\dots,c_n)$ be a feasible solution to the linear program, and define
\[
M = \sum_{i=0}^n c_i E_i.
\]

Observe that each $E_i$ admits a decomposition
\[
E_i = \varPhi_i^{} \varPhi_i^\dagger, 
\]
where the columns of $\varPhi_i$ are the vectors $\ket{\varphi_x}$ for $x \in S_i$, and the dagger \(\cdot^\dagger\) denotes the conjugate transpose operator. 
It follows that
\[
M = \sum_{i=0}^n c_i \varPhi_i^{} \varPhi_i^\dagger.
\]

Using the coefficients $c_i$, construct a matrix $N$ by concatenating $c_i$ copies of $\varPhi_i$ side by side for each $i$, namely,
\begin{equation*}
N = \Big(\underbrace{\varPhi_0, \dots, \varPhi_0}_{\text{$c_0$ copies}}, 
     \underbrace{\varPhi_1, \dots, \varPhi_1}_{\text{$c_1$ copies}}, 
     \dots, 
     \underbrace{\varPhi_n, \dots, \varPhi_n}_{\text{$c_n$ copies}}\Big).
\end{equation*}
Note that each $c_i$ is a natural number and that their sum $\sum_{i=0}^n c_i$ is positive, hence $N$ is not empty.
Let $k = \sum_{i=0}^n (q-1)^i \binom{n}{i} c_i$ be the number of columns of $N$, and define a map $\rho: \mathbb{G}^n \to \mathbb{C}^k$ 
that sends each vertex $x \in \mathbb{G}^n$ to the row of $N$ indexed by $x$.  
By construction, the complex inner product of $\rho(x)$ and $\rho(y)$ coincides with the $(x,y)$-entry of $M$, i.e.,
\[
\langle\rho(x),\rho(y)\rangle = M_{x,y}.
\]

On the other hand, applying the change-of-basis formula
\[
E_j = \frac{1}{q^n} \sum_{i=0}^n K_j(i) A_i,
\]
We can rewrite $M$ in terms of the adjacency matrices:
\[
M = \sum_{i=0}^n c_i E_i =\sum_{i=0}^{n}c_i\left(\frac{1}{q^n}\sum_{j=0}^nK_i(j)A_j\right)
  = \frac{1}{q^n} \sum_{j=0}^n \left( \sum_{i=0}^n c_i K_i(j) \right) A_j.
\]
By feasibility, the coefficient of $A_d$ in this expansion is zero. Hence, for any adjacent vertices $x,y \in H(n,q,d)$,
\[
\langle\rho(x),\rho(y)\rangle = M_{x,y} = 0,
\]
showing that $\rho$ is indeed an orthogonal representation of dimension $\sum_{i=0}^n (q-1)^i \binom{n}{i} c_i$, and is clearly modulus-one.
\end{proof}

\begin{corollary} \label{cor:OR-Case1}
If $d\geq\frac{(q-1)n}{q}$, then $\xi'(H(n,q,d))\leq qd$.
\end{corollary}

\begin{proof}
Since $d \ge \frac{(q-1)n}{q}$, we have $K_1(d) = (q-1)n - qd \le 0$.  
Therefore, by setting $c_0 = -K_1(d)$, $c_1 = 1$, and $c_i = 0$ for all other $i$, we obtain a feasible solution.  
Consequently, we have
\[
\xi'(H(n,q,d)) \le -K_1(d) + (q-1)\binom{n}{1} = qd - (q-1)n + (q-1)n = qd.
\]
\end{proof}

\begin{corollary}\label{cor:OR-Case2}
If $\frac{(q-1)n}{q}>d\geq\frac{(q-1)n-(\frac{q-2}{2})-\sqrt{(q-1)n+(\frac{q-2}{2})^2}}{q}$, then 
$$\xi'(H(n,q,d))\leq \left((q-1)n - \frac{q-2}{2}\right)qd -\frac{1}{2}q^2d^2.$$
\end{corollary}

\begin{proof}
Observe that
\begin{equation*}
    K_2(d) = \frac{1}{2}q^2d^2 - \left((q-1)n - \frac{q-2}{2}\right)qd+(q-1)^2\binom{n}{2},
\end{equation*}
which is negative if and only if
$\frac{(q-1)n-(\frac{q-2}{2})-\sqrt{(q-1)n+(\frac{q-2}{2})^2}}{q}< d <\frac{(q-1)n-(\frac{q-2}{2})+\sqrt{(q-1)n+(\frac{q-2}{2})^2}}{q}$.
Hence, when 
$\frac{(q-1)n}{q}>d\geq\frac{(q-1)n-(\frac{q-2}{2})-\sqrt{(q-1)n+(\frac{q-2}{2})^2}}{q}$, we have $K_2(d)\leq 0$. Therefore, setting $c_0 = -K_2(d)$, $c_2 = 1$, and $c_i = 0$ for all other $i$ yields a feasible solution.  

It follows that
\begin{equation*}
    \xi'(H(n,q,d)) \le -K_2(d) + (q-1)^2 \binom{n}{2}=\left((q-1)n - \frac{q-2}{2}\right)qd -\frac{1}{2}q^2d^2
\end{equation*}
\end{proof}

\begin{corollary}\label{cor:OR-Case3}
If $d=\delta n$ for some $0<\delta<\frac{q-1}{q}$, then for sufficiently large $n$, we have
\begin{equation*}
    \xi'(H(n,q,d))\leq q^{H_q\left(\frac{q-1-(q-2)\delta-2\sqrt{(q-1)\delta(1-\delta)}}{q}\right)n+o(n)}.
\end{equation*}
\end{corollary}

\begin{proof}
It is well known (see, e.g., \cite{levenshtein2002krawtchouk}) that the degree $d$ Krawtchouk polynomial $K_d(z)$ has $d$ distinct real zeros 
$z_1^{(d)} < \dots < z_d^{(d)}$, and that each interval between two consecutive zeros contains at least one integer. Moreover, $z_1^{(d)}$ is asymptotically given by
\begin{equation*}
        \lim_{n\to\infty}\frac{z_1^{(\delta n)}}{n} = \frac{q-1-(q-2)\delta-2\sqrt{(q-1)\delta(1-\delta)}}{q}.
\end{equation*}

Let $\lceil z_1^{(d)} \rceil$ denote the smallest integer greater than or equal to $z_1^{(d)}$.
From the properties listed above, we know that $\lceil z_1^{(d)} \rceil \in (z_1^{(d)}, z_2^{(d)})$.
Since $K_d(0) = (q - 1)^d \binom{n}{d} > 0$, by the definition of $\lceil z_1^{(d)} \rceil$ and the continuity of $K_d(z)$, it follows that $K_d(\lceil z_1^{(d)} \rceil) \le 0$. Hence, by the reciprocal law, we have $K_{\lceil z_1^{(d)} \rceil}(d) \le 0$.

Therefore, setting $c_0 = -K_{\lceil z_1^{(d)}\rceil}(d)$, $c_{\lceil z_1^{(d)}\rceil} = 1$, and $c_i = 0$ for all other $i$ yields a feasible solution.
Consequently, we have
\begin{equation*}
\begin{split}
\xi'(H(n,q,d))
&\leq -K_{\lceil z_1^{(d)}\rceil}(d) + (q-1)^{\lceil z_1^{(d)}\rceil} \binom{n}{\lceil z_1^{(d)}\rceil} \\
&\leq 2 (q-1)^{\lceil z_1^{(d)}\rceil} \binom{n}{\lceil z_1^{(d)}\rceil} \\
&\leq q^{H_q\left(\frac{q-1-(q-2)\delta-2\sqrt{(q-1)\delta(1-\delta)}}{q}\right)n +o(n)}, 
\end{split}
\end{equation*}
where the second inequality holds since $(q-1)^{\lceil z_1^{(d)}\rceil} \binom{n}{\lceil z_1^{(d)}\rceil}$ is the maximum eigenvalue in absolute value of $H(n,q,\lceil z_1^{(d)}\rceil)$ while $K_{\lceil z_1^{(d)}\rceil}(d)$ is another eigenvalue of this graph; the last inequality follows from the well-known entropy estimate for Hamming spheres, namely, $(q-1)^{t}\binom{n}{t}\leq q^{H_q(\frac{t}{n})n}$, where the estimate holds whenever $t\leq \frac{(q-1)n}{q}$.
\end{proof}

\begin{proof}[Proof of \cref{thm:qc-upper-Hamming}]
\cref{thm:qc-upper-Hamming} follows directly from \cref{lem:qc-upper} and the combination of \cref{cor:OR-Case1,cor:OR-Case2,cor:OR-Case3}.
\end{proof}

\subsection{Minimum eigenvalue of certain Hamming graphs}
\label{sec3:eigenvalue}
In this subsection, we derive a lower bound for the quantum chromatic number of $H(n,q,d)$ using \cref{lem:qc-lower}. Since the Hamming graph $H(n,q,d)$ is $(q-1)^d \binom{n}{d}$-regular, its maximum eigenvalue equals $(q-1)^d \binom{n}{d}$. Consequently, we only need to focus on the minimum eigenvalue of $H(n,q,d)$.

% The minimum eigenvalue is an important quantity relevant to many combinatorial problems, such as the maximum cut and intersecting families. However, determining it is generally a challenging task. Motivated by semidefinite programming approaches to the max-cut problem on Hamming graphs, Van Dam and Sotirov \cite{van2016new} conjectured that for $d \ge \frac{(q-1)n+1}{q}$, with $d$ even when $q=2$, the minimum eigenvalue of $H(n,q,d)$ is $K_d(1)$. Alon and Sudakov \cite{alon2000bipartite} proved this for $q=2$ with $n$ large and $d/n$ fixed. Dumer and Kapralova \cite{dumer2013spherically} proved it for $q=2$ and all $n$. Finally, this conjecture was proved by Brouwer et al. \cite{brouwer2018smallest}. We formulate it as follows:
For $d\geq\frac{(q-1)n+1}{q}$, Brouwer et al. \cite{brouwer2018smallest} showed that: 
\begin{lemma}[{\cite[Theorem 1.4]{brouwer2018smallest}}]\label{lem:eig-case1}
For $d \ge \frac{(q-1)n+1}{q}$ with $d$ even when $q=2$, the minimum eigenvalue of $H(n,q,d)$ is $K_d(1)$. %$K_d(1)=-(q-1)^d\binom{n}{d}\frac{qd-(q-1)n}{(q-1)n}$.
\end{lemma}

In the remaining part of this subsection, we focus on the minimum eigenvalue of $H(n,q,d)$ in the regime where $d$ lies slightly below the threshold $\frac{(q-1)n}{q}$, with $q$ fixed and $n$ sufficiently large. First, we have the following observation.  
\begin{lemma}
\label{lem:inequality_multiplicity}
Let $G$ be an $r$-regular graph. 
Let $\lambda_1, \lambda_2, \dots, \lambda_k$ be the eigenvalues of the adjacency matrix of $G$, with corresponding multiplicities $m_1, m_2, \dots, m_k$. 
Then, for each $i=1,2,\dots,k$, we have  
\[
\left(\frac{\lambda_i}{r}\right)^2 \le \frac{|V(G)|}{rm_i}.
\]
\end{lemma}

\begin{proof}
Let $A$ be the adjacency matrix of the $r$-regular graph $G$. 
Since each diagonal entry of $A^2$ equals $r$, we have $\mathrm{tr}(A^2)=|V(G)|r$. 
On the other hand, as the eigenvalues of $A^2$ are $\lambda_i^2$ with multiplicities $m_i$, we have $\mathrm{tr}(A^2)=\sum_{i=1}^k m_i\lambda_i^2$. 
Equating the two expressions gives $\sum_i m_i\lambda_i^2=|V(G)|r$. 
Hence $m_i\lambda_i^2\le |V(G)|r$ for each $i$, and dividing both sides by $m_i r^2$ yields $(\frac{\lambda_i}{r})^2 \le \frac{|V(G)|}{rm_i}$, as claimed.
\end{proof}

% \begin{lemma}\label{lem:eig-case2}
% For $\frac{(q-1)n}{q}\geq d\geq\frac{(q-1)n-\frac{q-2}{2}-\sqrt{(q-1)n+(\frac{q-2}{2})^2-(q-1)}+1}{q}$, with $d$ even when $q=2$ and $n$ sufficiently large, the minimum eigenvalue of $H(n,q,d)$ is $K_d(2)$.
% \end{lemma}

For the reader's convenience, here we restate the following lemma from the introduction.

\begin{llembis}{lem:eig-case2}
    For $\frac{(q-1)n}{q}\geq d\geq\frac{(q-1)n-\frac{q-2}{2}-\sqrt{(q-1)n+(\frac{q-2}{2})^2-(q-1)}+1}{q}$, with $d$ even when $q=2$ and $n$ sufficiently large, the minimum eigenvalue of $H(n,q,d)$ is $K_d(2)$.
\end{llembis}

\begin{proof}
As discussed in \cref{sec:pre}, the eigenvalues of $H(n,q,d)$ are given by $K_d(i)$ for $i=0,1,\dots,n$, with (formal) multiplicities $(q-1)^i\binom{n}{i}$. 

We begin by showing that $K_d(2)\leq K_d(i)$ for $i\leq 3$. By reciprocity,
$K_d(i)=(q-1)^d\binom{n}{d}\frac{K_i(d)}{(q-1)^i\binom{n}{i}}$,
and therefore the inequality $K_d(i)\leq K_d(j)$ is equivalent to
$\frac{K_i(d)}{(q-1)^i\binom{n}{i}}\leq \frac{K_j(d)}{(q-1)^j\binom{n}{j}}$.

Observe that 
\begin{equation*}
\begin{split}
K_0(d)=&\quad1, \\
K_1(d)=&\quad -qd+(q-1)n, \\
K_2(d)=&\quad\,\frac{1}{2}q^2d^2-\left((q-1)n-\frac{q-2}{2}\right)qd+(q-1)^2\binom{n}{2}, \\
K_3(d)=&-\frac{1}{6}q^3d^3+\frac{1}{2}\Big((q-1)n-(q-2)\Big) q^2d^{2}\\
    &-\frac{1}{6}\Big(3(q-1)^2n^2-3(2q-3)(q-1)n+2(q^2-3q+3)\Big) qd+(q-1)^3\binom{n}{3}.
\end{split}
\end{equation*}

Clearly, $K_0(d)>0$ for all $d$, $K_1(d)\geq0$ for $d\leq \frac{(q-1)n}{q}$, and $K_2(d)$ is negative if and only if $\frac{(q-1)n-\frac{q-2}{2}-\sqrt{(q-1)n+(\frac{q-2}{2})^2}}{q}<d< \frac{(q-1)n-\frac{q-2}{2}+\sqrt{(q-1)n+(\frac{q-2}{2})^2}}{q}$. Therefore, under the assumptions of the lemma on $d$, we have $K_2(d)<0$.

Next, we compare $\frac{K_2(d)}{(q-1)^2\binom{n}{2}}$ and $\frac{K_3(d)}{(q-1)^3\binom{n}{3}}$. Define
\begin{equation*}
f(z)=\frac{K_2(z)}{(q-1)^2\binom{n}{2}}-\frac{K_3(z)}{(q-1)^3\binom{n}{3}}.
\end{equation*}
Observe that $f(0)=1-1=0$, and hence the monomial $z$ divides $f(z)$. Consequently, $f(z)/z$ is a quadratic polynomial, whose roots are easily computed to be $\frac{(q-1)n-\frac{q-2}{2}\pm\sqrt{(q-1)n+(\frac{q-2}{2})^2-(q-2)}+1}{q}$. Since $f(z)$ is a cubic polynomial with positive leading coefficient, it follows that under the assumptions of the lemma on $d$, we have $f(d)\leq 0$. Thus, $K_d(2)\leq K_d(i)$ for $i\leq 3$.

It remains to show that $K_d(2)\leq K_d(i)$ for $4\leq i\leq n$. Note that when $q=2$ and $d$ is even, we have $K_d(i)=K_d(n-i)$. Therefore, in this case, it suffices to consider $4\leq i\leq\frac{n}{2}$.

In fact, we prove the following stronger statement:
\begin{itemize}
\item for $q=2$ and $4\leq i\leq\frac{n}{2}$;
\item for $q\geq 3$ and $4\leq i\leq n$,
\end{itemize}
we have $|K_d(2)|\geq |K_d(i)|$ for sufficiently large $n$.

By \cref{lem:inequality_multiplicity}, we have
\begin{equation*}
\left(\frac{K_d(i)}{(q-1)^{d}\binom{n}{d}}\right)^2 \leq \frac{q^n}{ (q-1)^{d}\binom{n}{d}(q-1)^i\binom{n}{i}}.
\end{equation*}
Thus, it suffices to show that
\begin{equation*}
\frac{q^n}{ (q-1)^{d}\binom{n}{d}(q-1)^i\binom{n}{i}} \le \left(\frac{K_d(2)}{(q-1)^d\binom{n}{d}}\right)^2.
\end{equation*}

On the one hand, we derive an upper bound for the left-hand side. Since $(q-1)^i\binom{n}{i}$ is unimodal in $i$ and attains its maximum at $i=\frac{(q-1)n}{q}$, we conclude that for sufficiently large $n$, $(q-1)^i\binom{n}{i}\geq (q-1)^4\binom{n}{4}$ for $q=2$ and $4\leq i\leq \frac{n}{2}$, and for $q\geq 3$ and $4\leq i\leq n$. Moreover, by a standard local limit estimate for the binomial distribution, for $d=\frac{(q-1)n}{q}-O(\sqrt{n})$, we have $(q-1)^d\binom{n}{d}\sim \frac{q^n}{\sqrt{n}}$. Therefore, under the assumptions of the lemma on $d$, we have
\begin{equation*}
\frac{q^n}{(q-1)^{d}\binom{n}{d}(q-1)^i\binom{n}{i}} \lesssim_q \frac{q^n}{\frac{q^n}{\sqrt{n}}(q-1)^4\binom{n}{4}} \lesssim_q n^{-3.5}.
\end{equation*}

On the other hand, under the assumptions of the lemma on $d$, the quadratic polynomial $K_2(d)$ attains its maximum at $d=\frac{(q-1)n-\frac{q-2}{2}-\sqrt{(q-1)n+(\frac{q-2}{2})^2-(q-1)}+1}{q}$, with value $-\frac{q-2}{2}-\sqrt{(q-1)n+(\frac{q-2}{2})^2-(q-1)}$. Thus,
\begin{equation*}
\left(\frac{K_d(2)}{(q-1)^d\binom{n}{d}}\right)^2 = \left(\frac{K_2(d)}{(q-1)^2\binom{n}{2}}\right)^2 \geq \left(\frac{-\frac{q-2}{2}-\sqrt{(q-1)n+(\frac{q-2}{2})^2-(q-1)}}{(q-1)^2\binom{n}{2}}\right)^2 \gtrsim_q n^{-3}.
\end{equation*}

Hence, we conclude that $K_d(2)\leq K_d(i)$ for all $i\geq 4$, which completes the proof.
\end{proof}

\begin{proof}[Proof of \cref{thm:qc-lower-Hamming}]
Together with \cref{lem:qc-lower,lem:eig-case1,lem:eig-case2}, this yields the proof.
\end{proof}

\section{Quantum chromatic number of generalized Hadamard graphs}\label{sec:Hadamard}
In this section, we discuss the quantum chromatic number of the generalized Hadamard graph $\varOmega_n^{(\mathbb{G})}$. These graphs admit a natural modules-one orthogonal representation of dimension $n$, which we briefly discuss in the following lemma.
\begin{lemma}\label{lem:qc-Haramard-upper}
$\xi'(\varOmega_n^{(\mathbb{G})})\leq n$.
\end{lemma}

\begin{proof}
Let $\varphi$ be a non-trival character of $\mathbb{G}$, consider the map $\rho:\mathbb{G}^n\to\mathbb{C}^n$ that send each vertex $x$ to the vector $(\varphi(x_1),\dots,\varphi(x_n))$, therefore for any two adjacent vertices $x,y$, i.e., $\comp(x-y)=\boldsymbol{n/q}$, we have
\begin{equation*}
\langle\rho(x),\rho(y)\rangle=\sum_{k=1}^n\overline{\varphi(x_k)}{\varphi(y_k)}=\sum_{k=1}\varphi(y_k-x_k)=\frac{n}{q}\sum_{g\in\mathbb{G}}\varphi(g)=0, 
\end{equation*}
where the last equality follows from the first orthogonality relation of characters.
It is clear that the orthogonal representation $\rho$ is modulus-one, since the values of $\varphi$ are $|\mathbb{G}|$-th roots of unity.
\end{proof}

In the remaining part of this section, we focus on the minimum eigenvalue of $\varOmega_n^{(\mathbb{G})}$, considering the cases $\mathbb{G} = \mathbb{Z}_q$ in \cref{sec4:cyclic group} and $\mathbb{G} = \mathbb{F}_q$ in \cref{sec4:finite field}.

\subsection{Minimum eigenvalue of \texorpdfstring{$\varOmega_n^{(\mathbb{Z}_q)}$}{}} \label{sec4:cyclic group}
Here, we apply the same approach as in \cref{sec3:eigenvalue} to determine the minimum eigenvalue of $\varOmega_n^{(\mathbb{Z}_q)}$.

We begin by establishing further properties of the generalized Krawtchouk polynomial $K^{(\mathbb{Z}_q)}_{\boldsymbol{i}}(\boldsymbol{j})$. %Throughout this subsection, we focus exclusively on $\mathbb{G} = \mathbb{Z}_q$ and, for simplicity, write $K_{\boldsymbol{i}}(\boldsymbol{j})$ instead of $K^{(\mathbb{Z}_q)}_{\boldsymbol{i}}(\boldsymbol{j})$.

From the discussion in \cref{sec:pre}, the eigenvalues of 
\(\varOmega_n^{(\mathbb{Z}_q)} = H_C(n, \mathbb{Z}_q, \boldsymbol{n/q})\) are given by
\(K^{(\mathbb{Z}_q)}_{\boldsymbol{n/q}}(\boldsymbol{r})\)
for all composition \(\boldsymbol{r} = (r_0, r_1, \dots, r_{q-1})\).

Let \(\boldsymbol{z} = (z_i : i \in \mathbb{Z}_q)\) be a tuple of indeterminates, and let
\[
C(\boldsymbol{z}) = (z_{j-i})_{i,j \in \mathbb{Z}_q} = 
\begin{pmatrix}
    z_0    & z_1 & \cdots & z_{q-2} & z_{q-1}\\ 
    z_{q-1}& z_0 & \cdots & z_{q-3} & z_{q-2}\\
    \vdots & \vdots & & \vdots & \vdots\\
    z_2 & z_3  & \cdots & z_0 & z_1 \\
    z_1 & z_2 & \cdots & z_{q-1} & z_0
\end{pmatrix}
\]
be the corresponding circulant matrix.

Observe that
\begin{equation*}
K^{(\mathbb{Z}_q)}_{\boldsymbol{r}}(\boldsymbol{n/q})
= [\boldsymbol{z}^{\boldsymbol{r}}]
\prod_{i\in \mathbb{Z}_q}\left(\sum_{j\in \mathbb{Z}_q} \zeta_q^{ij} z_j \right)^{n/q}
= [\boldsymbol{z}^{\boldsymbol{r}}]
\left(\prod_{i\in \mathbb{Z}_q} \sum_{j\in \mathbb{Z}_q} \zeta_q^{ij} z_j \right)^{n/q}
= [\boldsymbol{z}^{\boldsymbol{r}}]\Big(\det C(\boldsymbol{z})\Big)^{n/q},
\end{equation*}
where the last equality follows from the well-known property of circulant matrices: the determinant equals the product of its eigenvalues \(\sum_{j\in \mathbb{Z}_q} \zeta_q^{ij} z_j\) for \(i=0,1,\dots,q-1\).

Therefore, using the reciprocal property, we obtain
\begin{equation}\label{equ:det}
K^{(\mathbb{Z}_q)}_{\boldsymbol{n/q}}(\boldsymbol{r})
= \frac{\binom{n}{\boldsymbol{n/q}}}{\binom{n}{\boldsymbol{r}}} \, [\boldsymbol{z}^{\boldsymbol{r}}] \Big(\det C(\boldsymbol{z})\Big)^{n/q}.
\end{equation}

\begin{lemma}\label{lem:zero_eigenvalues}
    For any composition $\boldsymbol{r}=(r_0,\dots,r_{q-1})$, if 
    $\sum_{i=0}^{q-1} i r_i \not\equiv 0 \pmod q$, 
    then $K^{(\mathbb{Z}_q)}_{\boldsymbol{n/q}}(\boldsymbol{r}) = 0$.
\end{lemma}

\begin{proof}
 Let $e$ denote the all-one tuple in $\mathbb{Z}_q^n$, and let $a\in \mathbb{Z}_q^n$ be a tuple with $\comp(a) = \boldsymbol{r}$. Observe that 
 \[ e + S_{\boldsymbol{n/q}} = S_{\boldsymbol{n/q}}. \]
 Hence,
 \[ \zeta_q^{a\cdot e} \cdot \sum_{x \in S_{\boldsymbol{n/q}}} \zeta_q^{a \cdot x} 
    = \sum_{x \in S_{\boldsymbol{n/q}}} \zeta_q^{a \cdot (e+x)} 
    = \sum_{x \in S_{\boldsymbol{n/q}}} \zeta_q^{a \cdot x}. \]
 If $a \cdot e = \sum_{i \in \mathbb{Z}_q} i r_i \not\equiv 0 \pmod q$, then $\zeta_q^{a\cdot e} \neq 1$, which implies
 $K^{(\mathbb{Z}_q)}_{\boldsymbol{n/q}}(\boldsymbol{r}) = \sum_{x \in S_{\boldsymbol{n/q}}} \zeta_q^{a \cdot x} = 0$.
\end{proof}

\begin{lemma}\label{lem:symmetry}
Let $\boldsymbol{r} = (r_0, r_1, \dots, r_{q-1})$ be a composition, and let 
$\boldsymbol{r}^{(1)} = (r_1, r_2, \dots, r_{q-1}, r_0)$ denote its left cyclic shift. 
Then
\[
K^{(\mathbb{Z}_q)}_{\boldsymbol{n/q}}(\boldsymbol{r}^{(1)}) 
= (-1)^{\frac{(q-1)n}{q}} K^{(\mathbb{Z}_q)}_{\boldsymbol{n/q}}(\boldsymbol{r}).
\]
\end{lemma}

\begin{proof}
Let $\boldsymbol{z} = (z_i : i \in \mathbb{Z}_q)$ be a tuple of indeterminates, and let $\boldsymbol{z}^{(-1)} = (z_{q-1}, z_0, \dots, z_{q-2})$ denote its right cyclic shift. 
Then the circulant matrix $C(\boldsymbol{z})$ satisfies
\[
\det C(\boldsymbol{z}^{(-1)}) = (-1)^{q-1} \det C(\boldsymbol{z}).
\]
Observe that $\boldsymbol{z}^{\boldsymbol{r}^{(1)}}=z_0^{r_1}z_1^{r_2}\cdots z_{q-2}^{r_{q-1}}z_{q-1}^{r_0}=z_{q-1}^{r_0}z_0^{r_1}z_1^{r_2}\cdots z_{q-2}^{r_{q-1}}=(\boldsymbol{z}^{(-1)})^{\boldsymbol{r}}$, we obtain
\[
[\boldsymbol{z}^{\boldsymbol{r}^{(1)}}] \Big(\det C(\boldsymbol{z})\Big)^{n/q} 
= [(\boldsymbol{z}^{(-1)})^{\boldsymbol{r}}]\Big((-1)^{q-1} \det C(\boldsymbol{z}^{(-1)}) \Big)^{n/q} 
= (-1)^{\frac{(q-1)n}{q}} [\boldsymbol{z}^{\boldsymbol{r}}] \Big(\det C(\boldsymbol{z})\Big)^{n/q}.
\]
Therefore, 
\begin{equation*}
\begin{split}
K^{(\mathbb{Z}_q)}_{\boldsymbol{n/q}}(\boldsymbol{r}^{(1)}) 
&= \frac{\binom{n}{\boldsymbol{n/q}}}{\binom{n}{\boldsymbol{r}^{(1)}}} \, [\boldsymbol{z}^{\boldsymbol{r}^{(1)}}] \Big(\det C(\boldsymbol{z})\Big)^{n/q}\\
&=\frac{\binom{n}{\boldsymbol{n/q}}}{\binom{n}{\boldsymbol{r}}} \, (-1)^{\frac{(q-1)n}{q}}[\boldsymbol{z}^{\boldsymbol{r}}] \Big(\det C(\boldsymbol{z})\Big)^{n/q}=(-1)^{\frac{(q-1)n}{q}}K^{(\mathbb{Z}_q)}_{\boldsymbol{n/q}}(\boldsymbol{r}) .    
\end{split}
\end{equation*}
\end{proof}

\begin{lemma}\label{lem:eig-Zq}
Let $q\geq 2$ be a positive integer. Let $n$ be divisible by $q$ such that $\frac{(q-1)n}{q}$ is an even integer. 
Then, there exists a constant $N(q)$ such that if $n\ge N(q)$, then the minimum eigenvalue of $\varOmega_n^{(\mathbb{Z}_q)}$ is 
\[ K^{(\mathbb{Z}_q)}_{\boldsymbol{n/q}}(n-2,1,0,\dots,0,1) = -\frac{\binom{n}{n/q,\ldots,n/q}}{n-1}.
\]  
\end{lemma}

\begin{proof}
Let $\boldsymbol{r} = (r_0, r_1, \ldots, r_{q-1})$ be a composition.  
By Lemma \ref{lem:symmetry}, since $\frac{(q-1)n}{q}$ is even, the value of generalized Krawtchouk polynomial $K^{(\mathbb{Z}_q)}_{\boldsymbol{n/q}}$ evaluated at $\boldsymbol{r}$ is invariant under cyclic shifts of $\boldsymbol{r}$. 
Therefore, without loss of generality, we may assume that $r_0 \ge r_i$ for all $i = 1,2, \ldots, q-1$, which in particular implies $r_0 \ge n/q$.

\noindent\textbf{Case 1}. If $r_0 = n$, then $K^{(\mathbb{Z}_q)}_{\boldsymbol{n/q}}(n,0,\dots,0)=\binom{n}{n/q,\dots,n/q}$ is the maximum eigenvalue.

\noindent\textbf{Case 2}. If $r_0 = n-1$ and $r_k = 1$ for some $k\neq 0$, then $\sum_{i \in \mathbb{Z}_q }ir_i= k\not\equiv 0 \pmod q$, so by Lemma \ref{lem:zero_eigenvalues}, $K^{(\mathbb{Z}_q)}_{\boldsymbol{n/q}}(\boldsymbol{r}) = 0$.

\noindent\textbf{Case 3}. Suppose $r_0 = n-2$. Then, by Lemma \ref{lem:zero_eigenvalues}, the value $K^{(\mathbb{Z}_q)}_{\boldsymbol{n/q}}(\boldsymbol{r})$ is nonzero only if there exists some $k \in \mathbb{Z}_q \setminus \{0\}$ such that either $k \neq -k$ and $r_k = r_{-k} = 1$, or $k = -k$ and $r_k = 2$. Therefore, in both subcases, we have
\begin{align*}
    K^{(\mathbb{Z}_q)}_{\boldsymbol{n/q}}(\boldsymbol{r})&=\frac{\binom{n}{\boldsymbol{n/q}}}{\binom{n}{\boldsymbol{r}}}[z_0^{n-2}z_k z_{-k}]\Big(\det C(\boldsymbol{z})\Big)^{n/q}\\
&= \frac{\binom{n}{\boldsymbol{n/q}}}{\binom{n}{\boldsymbol{r}}}\binom{n/q}{1}[z_0^{q-2}z_kz_{-k}]\det C(\boldsymbol{z}) %= \frac{\binom{n}{n/q,\dots,n/q}}{q(n-1)}[z_0^{q-2}z_kz_{-k}]\det C(\boldsymbol{z}).
\end{align*}
Since $\det C(\boldsymbol{z})$ is the determinant of the circulant matrix $C(\boldsymbol{z})$, by definition, we have
\[\det C(\boldsymbol{z}) = \sum_{\sigma\in\mathfrak{S}_q}\prod_{i = 1}^qz_{\sigma(i)-i}, \]
where $\mathfrak{S}_q$ is the symmetric group of order $q$. 
The coefficient of $z_0^{q-2}$ in $\det C(\boldsymbol{z})$ only involve those $\sigma\in\mathfrak{S}_q$ that fix $q-2$ elements and  transpose two remaining elements. We have
\begin{align*}
    [z_0^{q-2}]\det C(\boldsymbol{z}) &= \sum_{\{k_1,k_2\}\subseteq\mathbb{Z}_q}(-1)z_{k_1-k_2}z_{k_2-k_1} \\
    &= \begin{cases}
    -q(z_1z_{-1}+z_{2}z_{-2}+\cdots+z_{\frac{q-1}{2}}z_{-\frac{q-1}{2}}), &\text{ if $q$ is odd,}\\
    -q(z_1z_{-1}+z_{2}z_{-2}+\cdots+z_{\frac{q}{2}-1}z_{-(\frac{q}{2}-1)})-\frac{q}{2}z_{\frac{q}{2}}^2, &\text{ otherwise}.
\end{cases} 
\end{align*}
Therefore, in both subcases, we obtain $K^{(\mathbb{Z}_q)}_{\boldsymbol{n/q}}(\boldsymbol{r})=-\frac{\binom{n}{n/q,\ldots,n/q}}{n-1}$.    

\noindent\textbf{Case 4}. Suppose $n/q\le r_0< n-\frac{q+3}{2}$. %It suffices to show that 
%\[\left|\frac{K^{(\mathbb{Z}_q)}_{\boldsymbol{n/q}}(\boldsymbol{r})}{\binom{n}{n/q,\ldots,n/q}}\right|\le \left|\frac{K^{(\mathbb{Z}_q)}_{\boldsymbol{n/q}}(n-2,1,0,\ldots,0,1)}{\binom{n}{n/q,\ldots,n/q}}\right| = \frac{1}{n-1}.\]
By Lemma~\ref{lem:inequality_multiplicity}, we have
\[\left(\frac{K^{(\mathbb{Z}_q)}_{\boldsymbol{n/q}}(\boldsymbol{r})}{\binom{n}{n/q,\ldots,n/q}}\right)^2 \le \frac{q^n}{\binom{n}{n/q,\ldots,n/q}\binom{n}{r_0,r_1,\cdots,r_{q-1}}}.\]
To show that $K^{(\mathbb{Z}_q)}_{\boldsymbol{n/q}}(n-2,1,0,\dots,0,1)$ is the minimum eigenvalue, it suffices to show that 
\[\frac{q^n}{\binom{n}{n/q,\ldots,n/q}\binom{n}{r_0,r_1,\dots,r_{q-1}}} \le \left(\frac{K^{(\mathbb{Z}_q)}_{\boldsymbol{n/q}}(n-2,1,0,\dots,0,1)}{\binom{n}{n/q,\ldots,n/q}}\right)^2,\]
since $K^{(\mathbb{Z}_q)}_{\boldsymbol{n/q}}(n-2,1,0,\dots,0,1)=-\frac{\binom{n}{n/q,\cdots,n/q}}{n-1}$, after ranging, it suffice to show that 
\begin{equation}
\label{eq:hadamard_general_cyclic_group_q}
    \binom{n}{r_0,r_1,\dots,r_{q-1}}\ge \frac{q^n(n-1)^2}{\binom{n}{n/q,\ldots,n/q}}.
\end{equation}
Since $\frac{q+3}{2}<n/q\le r_0< n-\frac{q+3}{2}$, 
\[\binom{n}{r_0,r_1,\ldots,r_{q-1}}\ge \binom{n}{r_0,n-r_0}\ge \binom{n}{\lfloor \frac{q+3}{2} \rfloor +1}\gtrsim_q n^{\lfloor\frac{q+3}{2} \rfloor +1}.\]
On the other hand, applying Stirling's approximation, we have
\[
\frac{q^n(n-1)^2}{\binom{n}{n/q,\ldots,n/q}} \le \frac{(n-1)^2\left(\sqrt{2\pi n/q}\right)^qe^{\frac{q^2}{12n}}}{\sqrt{2\pi n} }\lesssim_q n^{\frac{q+3}{2}}.
\]
Therefore, $\eqref{eq:hadamard_general_cyclic_group_q}$ is satisfied when $n$ is sufficiently large. 

\noindent\textbf{Case 5}. Suppose $n-\frac{q+3}{2}\le r_0\le n-3$. Let $3\le s\le \lfloor\frac{q+3}{2}\rfloor$ and $r_0 = n-s$. 
Let $\per C(\boldsymbol{z})$ denote the permanent of $C(\boldsymbol{z})$. Then, it is clear that
%Let $\det C(\boldsymbol{z})$ and $f_{per}$ denote the determinant and permanent of the circulant matrix with the first row as $[z_0\ z_1\ \cdots\ z_{q-1}]$, respectively. It is clear that
\begin{align}
\label{eq:big_proof_1}
    \Big|[z_0^{r_0}z_1^{r_1}\cdots z_{q-1}^{r_{q-1}}]\Big(\det C(\boldsymbol{z})\Big)^{n/q}\Big|\le [z_0^{r_0}z_1^{r_1}\cdots z_{q-1}^{r_{q-1}}]\Big(\per C(\boldsymbol{z})\Big)^{n/q}\le [z_0^{r_0}]\Big(\per C(z_0,1,\dots,1)\Big)^{n/q}.
\end{align}
    Since the coefficient of $z_0^{q-1}$ in the expansion of 
    \[\per C(z_0,1,\dots,1) = \per(\begin{bmatrix}
        z_0& 1&\cdots& 1\\
        1&z_0&\cdots&1\\
        \vdots&\vdots&\ddots &\vdots\\
        1& 1& \cdots& z_0
    \end{bmatrix})\]
    must be zero. Therefore, we can write 
    \[\per C(z_0,1,\ldots,1) = a_qz_0^q+a_{q-1}z_0^{q-1}+a_{q-2}z_0^{q-2}+\cdots+a_0, \]
    for some nonnegative integers $a_{q} = 1,a_{q-1} = 0 $, and $ a_i\le q!$, $i = 0,1,\ldots,q-2$. Then
    \begin{align}
    \label{eq:big_proof_2}
        [z_0^{r_0}]\Big(\per C(z_0,1,\ldots,1)\Big)^{n/q} = [z_0^{n-s}]\left(\sum_{i = 0}^qa_iz_0^i\right)^{n/q} = \sum_{\substack{k_1+k_2+\cdots+k_{n/q} =n-s\\ \text{each } k_j\neq q-1}}a_{k_1}a_{k_2}\cdots a_{k_{n/q}}.
    \end{align}
    Let $\boldsymbol{k} = (k_1,k_2,\ldots,k_{n/q})$ denote a feasible solution to $\sum_{j=1}^{n/q}k_j = n-s $ with $ k_j\neq q-1$ for all $j\in [n/q]$, and let $N(\boldsymbol{k})$ be the number of $j\in [n/q]$ such that $k_j\le q-2$. Since $k_j\neq q-1$ for all $j\in [n/q]$, we have $n-s = \sum_{j = 1}^{n/q}k_j\le q-2N(\boldsymbol{k})$, implying that $N(\boldsymbol{k})\le \frac{s}{2}$. Therefore,
    \begin{equation}
    \begin{aligned}
        \label{eq:big_proof_3}\sum_{\substack{k_1+k_2+\cdots+k_{n/q} =n-i\\ \text{each }k_j\neq q-1}}a_{k_1}a_{k_2}\cdots a_{k_{n/q}} &= \sum_{\ell = 0}^{\lfloor \frac{s}{2}\rfloor}\sum_{N(\boldsymbol{k}) = \ell}a_{k_1}a_{k_2}\cdots a_{k_{n/q}}\\
    &\le \sum_{\ell = 0}^{\lfloor \frac{s}{2}\rfloor}\sum_{N(\boldsymbol{k}) = \ell}(q!)^\ell\le \sum_{\ell = 0}^{\lfloor \frac{s}{2}\rfloor}\binom{n/q}{\ell}(q-1)^\ell(q!)^\ell\lesssim_q n^{\frac{s}{2}}.
    \end{aligned}
    \end{equation}
    Combining \eqref{eq:big_proof_1},\eqref{eq:big_proof_2} and \eqref{eq:big_proof_3}, we obtain
    \[\Big|[z_0^{r_0}z_1^{r_1}\cdots z_{q-1}^{r_{q-1}}]\Big(\det C(\boldsymbol{z})\Big)^{n/q}\Big|\lesssim_q n^{\frac{s}{2}}.\]
    Therefore,
    \[ \left|\frac{K^{(\mathbb{Z}_q)}_{\boldsymbol{n/q}}(\boldsymbol{r})}{\binom{n}{n/q,\dots,n/q}}\right| = \frac{\Big|[z_0^{r_0}z_1^{r_1}\cdots z_{q-1}^{r_{q-1}}]\Big(\det C(\boldsymbol{z})\Big)^{n/q}\Big|}{\binom{n}{n-s,r_1,\dots,r_{q-1}}}\lesssim_q\frac{n^{\frac{s}{2}}}{\binom{n}{s}}\lesssim_q n^{-\frac{s}{2}}, \]
    which decreases faster than
    \begin{equation*}
        \left|\frac{K^{(\mathbb{Z}_q)}_{\boldsymbol{n/q}}(n-2,1,0,\dots,0,1)}{\binom{n}{n/q,\dots,n/q}}\right|=\frac{1}{n-1}. 
    \end{equation*}
    Hence $\left|K^{(\mathbb{Z}_q)}_{\boldsymbol{n/q}}(\boldsymbol{r})\right|\le \left|K^{(\mathbb{Z}_q)}_{\boldsymbol{n/q}}(n-2,1,0,\dots,0,1)\right|$ for sufficiently large $n$.
\end{proof}

\subsection{Minimum eigenvalue of \texorpdfstring{$\varOmega_n^{(\mathbb{F}_q)}$}{}} \label{sec4:finite field}
Here, we present an argument of algebraic flavor to determine the minimum eigenvalue of $\varOmega_n^{(\mathbb{F}_q)}$. We first consider the case $n = q$, and then lift it to the prime power case where $n$ is divisible by $q$.

\begin{lemma}\label{lem:eig-Fq^q}
For any prime power $q$, the minimum eigenvalue of $\varOmega_q^{(\mathbb{F}_q)}$ is $-\frac{q!}{q-1}$.
\end{lemma}

\begin{proof}
Recall that $\varOmega_q^{(\mathbb{F}_q)} = \mathrm{Cay}(\mathbb{F}_q^q, S_{\boldsymbol{1}})$,
where $S_{\boldsymbol{1}}$ is the set of vectors in $\mathbb{F}_q^q$ having composition $\boldsymbol{1} = (1,1,\ldots,1)$.

Suppose that $q$ is a power of a prime $p$. By Theorem~\ref{thm:eigcay}, the eigenvalues of $\varOmega_q^{(\mathbb{F}_q)}$ are
\begin{equation*}
    \mu_a=\sum_{s\in S_{\boldsymbol{1}}}\zeta_p^{\operatorname{Tr}(s\cdot a)}, \quad a\in\mathbb{F}_q^q.
\end{equation*}
Observe that for each $x\in\mathbb{F}_q^\times$, we have $xS_{\boldsymbol{1}}=\{xs:s\in S_{\boldsymbol{1}}\}=S_{\boldsymbol{1}}$. Therefore, 
\begin{equation*}
\begin{split}
    \mu_a=\sum_{s\in S_{\boldsymbol{1}}}\zeta_p^{\operatorname{Tr}(a\cdot s)}
    &=\frac{1}{q-1}\sum_{x\in\mathbb{F}_q^\times}\sum_{s\in S_{\boldsymbol{1}}}\zeta_p^{\operatorname{Tr}(a\cdot (xs))}\\
    &=\frac{1}{q-1}\sum_{s\in S_{\boldsymbol{1}}}\sum_{x\in\mathbb{F}_q^\times}\zeta_p^{\operatorname{Tr}(x(a\cdot s))}\\
    &=\frac{1}{q-1}\left(\sum_{s\in S_{\boldsymbol{1}}:a\cdot s=0}(q-1)+\sum_{s\in S_{\boldsymbol{1}}:a\cdot s\neq0}(-1)\right)\\
    &=\frac{1}{q-1}\left(\sum_{s\in S_{\boldsymbol{1}}:a\cdot s=0}q-\sum_{s\in S_{\boldsymbol{1}}}1\right)\ge -\frac{|S_{\boldsymbol{1}}|}{q-1}, 
\end{split}
\end{equation*}
where the last inequality becomes equality if $a = (1,-1,0,\ldots,0)$, since there is no $s\in S$ that satisfies $a\cdot s = 0$. 
Thus, the minimum eigenvalue of $\varOmega_q^{(\mathbb{F}_q)}$ is $ -\frac{q!}{q-1}$.
\end{proof}

\begin{lemma}
\label{lem:eig-Fq^n}
Let $q$ and $n$ be prime powers with $q\mid n$.
Then, the minimum eigenvalue of $\varOmega_n^{(\mathbb{F}_q)}$ is $-\frac{\binom{n}{n/q,\dots,n/q}}{n-1}$.
\end{lemma}

\begin{proof}
Recall that $\varOmega_n^{(\mathbb{F}_n)} = \mathrm{Cay}(\mathbb{F}_n^n, S_{\boldsymbol{1}})$ and $\varOmega_n^{(\mathbb{F}_q)} = \mathrm{Cay}(\mathbb{F}_q^n, S_{\boldsymbol{n/q}})$.

Let $\theta_0:\mathbb{F}_n\rightarrow \mathbb{F}_q$ be any surjective linear map over the prime field, and let $\theta:\mathbb{F}_n^n\rightarrow \mathbb{F}_q^n$ be the map induced by $\theta_0$, defined componentwise by
\[
\theta(a_1,\dots,a_n) = (\theta_0(a_1),\dots,\theta_0(a_n)).
\]

For a character $\psi$ of $\mathbb{F}_q^n$, we write $\psi(S) := \sum_{s\in S} \psi(s)$ for brevity. 
Since every element of $S_{\boldsymbol{n/q}}$ has the same number of preimages in $S_{\boldsymbol{1}}$, it follows that
\[
(\psi \circ \theta)(S_{\boldsymbol{1}}) = \frac{|S_{\boldsymbol{1}}|}{|S_{\boldsymbol{n/q}}|} \, \psi(S_{\boldsymbol{n/q}}).
\]

Since $\psi \circ \theta$ is a character of $\mathbb{F}_n^n$, $(\psi \circ \theta)(S_{\boldsymbol{1}})$ is an eigenvalue of $\varOmega_n^{(\mathbb{F}_n)}$. By \cref{lem:eig-Fq^q}, we have
\[
(\psi \circ \theta)(S_{\boldsymbol{1}}) \ge -\frac{|S_{\boldsymbol{1}}|}{n-1},
\]
which immediately implies
\[
\psi(S_{\boldsymbol{n/q}}) \ge -\frac{|S_{\boldsymbol{n/q}}|}{n-1}.
\]
Thus, every eigenvalue of $\varOmega_n^{(\mathbb{F}_q)}$ is lower bounded by $-\frac{|S_{\boldsymbol{n/q}}|}{n-1}$.

To see that this bound is attained, consider $a = (1,-1,0,\dots,0)\in \mathbb{F}_q^n$ and the character $\psi_a$ defined by $\psi_a(x) = \zeta_p^{\mathrm{Tr}(a\cdot x)}$.  

Observe that for each $x \in \mathbb{F}_q^\times$, $xS_{\boldsymbol{n/q}} = S_{\boldsymbol{n/q}}$. In other words, there is a natural action of the multiplicative group $\mathbb{F}_q^\times$ on $S_{\boldsymbol{n/q}}$ by scalar multiplication. Denote the set of orbits by $S_{\boldsymbol{n/q}}/\mathbb{F}_q^\times$; each orbit contains exactly $q-1$ elements, and these orbits partition $S_{\boldsymbol{n/q}}$.  

Then we can write
\begin{align*}
\psi_a(S_{\boldsymbol{n/q}})
&= \sum_{x \in S_{\boldsymbol{n/q}}} \zeta_p^{\mathrm{Tr}(x_1 - x_2)} \\
&= \sum_{[x] \in S_{\boldsymbol{n/q}}/\mathbb{F}_q^\times} \sum_{\alpha \in \mathbb{F}_q^\times} \zeta_p^{\mathrm{Tr}(\alpha(x_1 - x_2))} \\
&= \sum_{[x] \in S_{\boldsymbol{n/q}}/\mathbb{F}_q^\times: x_1 = x_2} (q-1) + \sum_{[x] \in S_{\boldsymbol{n/q}}/\mathbb{F}_q^\times: x_1 \neq x_2} (-1) \\
&= \sum_{[x] \in S_{\boldsymbol{n/q}}/\mathbb{F}_q^\times: x_1 = x_2} q - \sum_{[x]\in S_{\boldsymbol{n/q}}/\mathbb{F}_q^\times} 1 \\
&= \frac{q-1}{q} \bigl|\{ x \in S_{\boldsymbol{n/q}} : x_1 = x_2 \}\bigr| - \frac{|S_{\boldsymbol{n/q}}|}{q-1} \\
&= \frac{q-1}{q} \binom{q}{1} \binom{n-2}{n/q-2, n/q, \dots, n/q} - \frac{|S_{\boldsymbol{n/q}}|}{q-1} \\
&= -\frac{|S_{\boldsymbol{n/q}}|}{n-1}.
\end{align*}
Hence, the lower bound is attained, completing the proof.
\end{proof}

\begin{proof}[Proof of \cref{thm:qc-Hadamard}]
\cref{thm:qc-Hadamard} follows from \cref{lem:qc-upper,lem:qc-Haramard-upper,lem:qc-lower,lem:eig-Zq,lem:eig-Fq^n}.
\end{proof}

\section{Chromatic number of generalized Hadamard graphs}\label{sec:Hadamard-sep}

\subsection{Chromatic number of \texorpdfstring{$\varOmega_n^{(\mathbb{Z}_q)}$}{}}

In this subsection, we determine the independence number of the generalized Hadamard graph $\varOmega_n^{(\mathbb{Z}_q)}$. Our approach relies on the seminal work of Frankl and R\"odl \cite{frankl1987forbidden}, who utilized the method of \emph{forbidden intersection patterns} to bound the independence numbers of Hamming graphs. We adapt their machinery to the context of generalized Hadamard graphs.

First, we recall the necessary definitions. For positive integers $k_1,\dots,k_s$ summing to $n$, let $\binom{[n]}{k_1,\dots,k_s}$ denote the set of ordered partitions $\boldsymbol{A}=(A_1,\dots,A_s)$ of $[n]$ where $|A_i|=k_i$. The \emph{intersection pattern} of two partitions $\boldsymbol{A},\boldsymbol{B} \in \binom{[n]}{k_1,\dots,k_s}$ is the $s \times s$ matrix $M=(m_{ij})_{i\in [s],j\in [s]}$ defined by $m_{ij}=|A_i\cap B_j|$. A matrix $M$ is called \textit{legitimate} with respect to the types $(k_0,\ldots,k_s)$ and $(\ell_0,\ldots,\ell_s)$ if its row and column sums satisfy: for each $i\in [s]$, $\sum_{j\in [s]} m_{ij}=k_i$; for each $j\in [s]$, $\sum_{i\in [s]} m_{ij}=\ell_j$.

% First, we recall the necessary definitions. For positive integers $k_1,\dots,k_s$ summing to $n$, let $\binom{[n]}{k_1,\dots,k_s}$ denote the set of ordered partitions $\boldsymbol{A}=(A_1,\dots,A_s)$ of $[n]$ where $|A_i|=k_i$. The \emph{intersection pattern} of two partitions $\boldsymbol{A} \in \binom{[n]}{k_1,\dots,k_s}$ and $\boldsymbol{B} \in \binom{[n]}{\ell_1,\dots,\ell_t}$ is the $s \times t$ matrix $M=(m_{ij})_{i\in [s],j\in [t]}$ defined by $m_{ij}=|A_i\cap B_j|$.
% A matrix $M$ is called \textit{legitimate} with respect to the types $(k_0,\ldots,k_s)$ and $(\ell_0,\ldots,\ell_t)$ if its row and column sums satisfy: for each $i\in [s]$, $\sum_{j\in [t]} m_{ij}=k_i$; for each $j\in [t]$, $\sum_{i\in [s]} m_{ij}=\ell_j$.

\begin{theorem}[{\cite[Theorem 1.15]{frankl1987forbidden}}]\label{thm:frankl1987forbidden}
Let $\eta$ be an arbitrary positive constant, 
$\mathcal{F}\subseteq\binom{[n]}{k_1,\dots,k_s}$, and let $M=(m_{ij})$ be a legitimate intersection pattern satisfying $m_{ij} \geq \eta n$.
Then there exists a positive constant $\varepsilon=\varepsilon(\eta)$ such that whenever
\begin{equation*}
    |\mathcal{F}|>(1-\varepsilon)^n\binom{n}{k_1,\dots,k_s}, 
\end{equation*}
$\mathcal{F}$ contains distinct $\boldsymbol{F}, \boldsymbol{F}'\in\mathcal{F}$ satisfying
\begin{equation*}
    |\boldsymbol{F}\cap \boldsymbol{F}'|=M.
\end{equation*}
\end{theorem}

For a vector $x \in \mathbb{Z}_q^n$, we associate the partition $\boldsymbol{A}(x) = (A_0(x), \dots, A_{q-1}(x))$ where $A_t(x) = \{ i \in [n] : x_i = t \}$ for $t\in \mathbb{Z}_q$. The \emph{composition} of $x$ is the tuple $\mathrm{comp}(x) = (|A_0(x)|, \dots, |A_{q-1}(x)|)$. Crucially, for any $x, y \in \mathbb{Z}_q^n$, the composition of their difference is determined by the intersection pattern $M = |\boldsymbol{A}(x) \cap \boldsymbol{A}(y)|$ via the relation
\begin{align}
\label{eq:relate_composition_and_pattern}
    |A_t(x-y)| = \sum_{j\in \mathbb{Z}_q} m_{(j+t)j},\quad t\in \mathbb{Z}_q.
\end{align}
To apply Theorem \ref{thm:frankl1987forbidden}, we must construct legitimate intersection patterns that force a specific difference composition (namely, the constant tuple $(n/q, \dots, n/q)$, which corresponds to adjacency in $\varOmega_n^{(\mathbb{Z}_q)}$). This construction is guaranteed by the following lemma.

\begin{lemma}
\label{lem:construct_intersection_pattern}
    Let $\mathbb{A}$ be the set of $q \times q$ matrices over $\mathbb{Z}_q$ defined by
    \[
        \mathbb{A} = \left\{ A = (x_i + x_j)_{i,j \in \mathbb{Z}_q} : x_i \in \mathbb{Z}_q \text{ and } \sum_{i \in \mathbb{Z}_q} x_i = 0 \right\}.
    \]
    For any $A = (a_{ij})_{i,j\in \mathbb{Z}_q} \in \mathbb{A}$, there exists a $q \times q$ integer matrix $L = (l_{ij})_{i,j\in \mathbb{Z}_q}$ satisfying:
    \begin{enumerate}
        \item[(1)] $l_{ij} \equiv a_{ij} \pmod{q}$ for all $i, j \in \mathbb{Z}_q$;
        \item[(2)] $\sum_{j\in \mathbb{Z}_q} l_{ij} = 0$ for each $i\in \mathbb{Z}_q$; and $\sum_{i\in \mathbb{Z}_q} l_{ij} = 0$ for each $j\in \mathbb{Z}_q$;
        \item[(3)] $\sum_{j\in \mathbb{Z}_q} l_{(j+t)j} = 0$ for each $t \in \mathbb{Z}_q$;
        \item[(4)] $|l_{ij}| < q^3$ for all $i, j \in \mathbb{Z}_q$.
    \end{enumerate}
\end{lemma}

\begin{proof}
    Consider the code $C = \{ x = (x_0,\ldots,x_{q-1}) \in \mathbb{Z}_q^q : \sum_{i\in \mathbb{Z}_q} x_i = 0 \}$. It is a free $\mathbb{Z}_q$-module spanned by the basis $v_0 = [1, -1, 0, \dots, 0]$ and its cyclic shifts, denoted by $v_s$ ($s=1, \dots, q-2$). The map $\tau: C \to \mathbb{A}$ given by $x \mapsto (x_i + x_j)_{i,j\in \mathbb{Z}_q}$ is $\mathbb{Z}_q$-linear and surjective.
    
    We first lift $\tau(v_0)$. Define the integer matrix $L_0$ by
    \[
        L_0 = \begin{pmatrix}
            2-q & 0 & 1 & \cdots & 1 \\
            0 & q-2 & -1 & \cdots & -1 \\
            1 & -1 & 0 & \cdots & 0 \\
            \vdots & \vdots & \vdots & \ddots & \vdots \\
            1 & -1 & 0 & \cdots & 0 
        \end{pmatrix}.
    \]
    Observe that $L_0 \equiv \tau(v_0) \pmod q$. One can verify that $L_0$ satisfies conditions (2)--(4), and the absolute values of its entries are bounded by $q$. For $1\le s \le q-2$, let $L_s$ be the matrix obtained by cyclically shifting the rows and columns of $L_0$ $s$ times respectively. By translation invariance, each $L_s$ is a valid lift for $\tau(v_s)$.

    Now, for arbitrary $A \in \mathbb{A}$, we can write $A = \sum_{s = 0}^{q-2} c_s \tau(v_s)$ for some $c_s\in \mathbb{Z}_q$. Let $C_s \in \{0, \dots, q-1\}$ be the integer lift of $c_s$. The matrix $L = \sum_{s_0}^{q-2} C_s L_s$ satisfies $L \equiv A \pmod q$. Since conditions (2)--(3) are linear and homogeneous, they are preserved by linear combinations. Finally, the entries are bounded by $|l_{ij}| \le \sum_{s = 0}^{q-2} |C_s| |(L_s)_{ij}| \le (q-1)q^2 < q^3$.
\end{proof}

\begin{lemma}\label{lem:Zq ind-number}
    Let $q \ge 2$ be an integer. There exist constants $\delta(q) > 0$ and $N(q)$ such that for all $n \geq N(q)$ with $\frac{(q-1)n}{q}$ being an even integer,
    \begin{equation*}
    \alpha(\varOmega_n^{(\mathbb{Z}_q)}) \leq (q-\delta)^n.
    \end{equation*}
\end{lemma}

\begin{proof}
    Let $C$ be an independent set of $\varOmega_n^{(\mathbb{Z}_q)}$. We partition $C$ by composition:
    \[
        C = \bigcup_{k_0+\cdots+k_{q-1} = n} C^{(k_0,\ldots,k_{q-1})}, \quad \text{where } C^{(k_0,\ldots,k_{q-1})} = \{x \in C : \mathrm{comp}(x)= (k_0, \dots, k_{q-1})\}.
    \]
    Set $\beta = \frac{2}{3q}$. We analyze two cases based on the balance of the composition.

   \textbf{Case 1.} Suppose $k := \min \{k_0, \ldots, k_{q-1}\} < \beta n$. Recall that the $q$-ary entropy function $H_q(x)= x\log_q(q-1)-x\log_q(x)-(1-x)\log_q(1-x)$ is strictly decreasing on the interval $[\frac{q-1}{q}, 1]$. Since $\beta < 1/q$, it follows that $\frac{n-k}{n} > 1-\beta > \frac{q-1}{q}$, which implies $H_q(1-\beta) < H_q(\frac{q-1}{q}) = 1$. Using the standard estimate for the binomial coefficient, we obtain
    \begin{align}
    \label{eq:unbalanced_upper_bound}
    |C^{(k_0,\ldots,k_{q-1})}| \le \binom{n}{n-k}(q-1)^{n-k} \le 2^{nH_2(\frac{n-k}{n})}(q-1)^{n-k}= q^{n H_q\left(\frac{n-k}{n}\right)} \le q^{n H_q(1-\beta)}.
    \end{align}
    Since $H_q(1-\beta) < 1$, this upper bound is exponentially small relative to $q^n$.

    \textbf{Case 2.} Suppose $k_i \ge \beta n$ for all $i = 0,\ldots,q-1$.
    We construct a legitimate forbidden intersection pattern $M$ with respect to $(k_0,\ldots,k_{q-1})$ and then apply Theorem \ref{thm:frankl1987forbidden}. Let $A = (a_{ij})_{i,j\in \mathbb{Z}_q}$ be a $q\times q$ matrix over $\mathbb{Z}_q$ with entries $a_{ij} = k_i + k_j - n/q\mod q$, where an integer modulo $q$, say $y \mod q$, represents the coset $y+\langle q\rangle\in \mathbb{Z}_q$.
    
    First, we verify $A \in \mathbb{A}$ (defined in Lemma \ref{lem:construct_intersection_pattern}), which is equivalent to verifying the existence of $x_0, \ldots, x_{q-1} \in \mathbb{Z}_q$ such that $\sum x_i = 0$ and $A = (x_i + x_j)_{i,j \in \mathbb{Z}_q}$.
       \begin{itemize}
        \item If $q$ is odd, 2 is invertible. Note that $a_{ij} \equiv k_i+k_j-\frac{n}{q} \equiv (k_i-\frac{q+1}{2}\cdot \frac{n}{q})+(k_j-\frac{q+1}{2}\cdot \frac{n}{q})\pmod{q}$. Let $x_i = k_i - \frac{q+1}{2}\cdot \frac{n}{q} \mod{q}$, for $i\in \mathbb{Z}_q$. The sum is $\sum_{i\in \mathbb{Z}_q} x_i\equiv \sum_{i\in \mathbb{Z}_q} k_i \equiv n \equiv 0 \pmod{q}$. 
        \item If $q$ is even, the condition  that $(q-1)n/q$ is even implies $n/q$ is even. Note that $a_{ij}\equiv k_i+k_j-\frac{n}{q} \equiv (k_i-\frac{n}{2q})+(k_j-\frac{n}{2q})$. Let  $x_i = k_i - \frac{q}{2n}$, for $i\in \mathbb{Z}_q$. The sum is $\sum_{i\in \mathbb{Z}_q} x_i \equiv \sum_{i\in \mathbb{Z}_q} k_i \equiv 0 \pmod{q}$. 
    \end{itemize}
    In either case, we have $A = (x_i+x_j)_{i,j\in \mathbb{Z}_q}$ with $\sum_{i\in \mathbb{Z}_q} x_i= 0$, and thus $A\in \mathbb{A}$.
    
    By Lemma \ref{lem:construct_intersection_pattern}, there is an integer matrix $L = (l_{ij})_{i,j\in \mathbb{Z}_q}$ having the properties: $L\equiv A\pmod{q}$; the row sums, column sums and diagonal sums of $L$ are zero; every entry $l_{ij}$ satisfies $|l_{ij}|\le q^3$.  Define $M = (m_{ij})_{i,j\in \mathbb{Z}_q}$ by
    \[
        m_{ij} = \frac{k_i + k_j - n/q - l_{ij}}{q}.
    \]
    Since $A-L \equiv \mathbf{0} \pmod q$, $m_{ij}$ are integers. We verify the following properties of $M$:
    \begin{enumerate}
        \item[(i)] For each $i\in \mathbb{Z}_q$, we have 
        \[\sum_{j\in \mathbb{Z}_q} m_{ij} = \frac{1}{q}\left(q k_i + \sum_{j\in \mathbb{Z}_q} k_j - n - \sum_{j\in \mathbb{Z}_q} l_{ij}\right) = \frac{1}{q}(q k_i + n - n - 0) = k_i.\]
        Symmetry implies that for each $j\in \mathbb{Z}_q$, we have $\sum_{i\in \mathbb{Z}_q}m_{ij} = k_j$.
        \item[(ii)] For each $t \in \mathbb{Z}_q$, we have
        \[
            \sum_{j\in \mathbb{Z}_q} m_{(j+t)j} = \frac{1}{q} \left( \sum_{j\in \mathbb{Z}_q} (k_{j+t}+k_j) - \sum_{j\in \mathbb{Z}_q} \frac{n}{q} - \sum_{j\in \mathbb{Z}_q} l_{(j+t)j} \right) = \frac{1}{q}(2n - n - 0) = \frac{n}{q}.
        \]
    \end{enumerate}
    Property (i) means that $M$ is a legitimate intersection pattern with respect to $(k_0,\ldots,k_{q-1})$. By the relation \eqref{eq:relate_composition_and_pattern}, (ii) implies that any two vectors $x, y$ with $|\boldsymbol{A}(x)\cap \boldsymbol{A}(y)| = M$ satisfy $\mathrm{comp}(x-y) = (n/q, \dots, n/q)$. Since $C$ is an independent set, no two vectors in $C$ can have this difference; thus, $M$ is a forbidden intersection pattern for $C^{(k_0,\ldots,k_{q-1})}$.
    
    Finally, using $|l_{ij}| < q^3$ and $k_i \ge \beta n$, for $n\ge 12q^4$, we have
    \[
        m_{ij} \ge \frac{2\beta n}{q} - \frac{n}{q^2} - q^2=\frac{4n}{3q^2}-\frac{n}{q^2}-q^2 \ge \frac{n}{4q^2}.
    \]
    It follows from Theorem \ref{thm:frankl1987forbidden} (with $\eta = 1/(4q^2)$ and $s = q$) that there is some $\varepsilon = \varepsilon(q) > 0$ satisfying
    \begin{align}
    \label{eq:balanced_upper_bound}
        |C^{(k_0,\ldots,k_{q-1})}| \le (1-\varepsilon)^n \binom{n}{k_0,\ldots,k_{q-1}}\le (1-\varepsilon)^n q^n.
    \end{align}

   Combining the upper bounds \eqref{eq:unbalanced_upper_bound} and \eqref{eq:balanced_upper_bound} from the two cases above, we obtain that for $n \ge 12q^4$,$$|C| \le \sum_{k_0+\cdots+k_{q-1} = n} |C^{(k_0,\ldots,k_{q-1})}| \le \binom{n+q-1}{q-1} \left( \max\{q^{H_q(1-\beta)}, q(1-\varepsilon)\} \right)^n. $$Consequently, there exist constants $\delta = \delta(q) > 0$ and $N = N(q)$ such that for all $n \ge N$ where $(q-1)n/q$ is an even integer, we have$$|C| \le (q-\delta)^n.$$This completes the proof.
\end{proof}

\subsection{Chromatic number of \texorpdfstring{$\varOmega_n^{(\mathbb{F}_q)}$}{}}
In this subsection, we establish a homomorphism from the Hamming graph to the generalized Hadamard graph $\varOmega_n^{(\mathbb{F}_q)}$.

\begin{lemma}\label{lem: homo hamming to Fq}
Let $q$ be a prime power and $n, d$ be positive integers such that $n = qd$. Let $m = \lfloor n/(q-1) \rfloor$. Then there exists a graph homomorphism from the Hamming graph $H(m, q, d)$ to $\varOmega^{(\mathbb{F}_q)}_{n}$.
\end{lemma}

\begin{proof}
We treat the Hamming graph $H(m, q, d)$ as the Cayley graph $\cay(\mathbb{F}_q^m,S_d)$, where $S_d$ consists of the tuple whose Hamming weight is $d$. Let $\mathbb{F}_q^\times = \{\lambda_1, \lambda_2, \dots, \lambda_{q-1}\}$ denote the set of all non-zero elements in the finite field $\mathbb{F}_q$.  Let $M \in \mathbb{F}_q^{m \times n}$ be a block matrix  defined by
\[
    M = \begin{bmatrix}
        \lambda_1 {I}_m & \lambda_2 {I}_m & \cdots & \lambda_{q-1} {I}_m & {0}_{m \times (n - (q-1)m)}
    \end{bmatrix},
\]
where ${I}_m$ is the $m \times m$ identity matrix and ${0}$ is the zero matrix of the appropriate dimensions. We define a mapping $\psi: \mathbb{F}_q^m \to \mathbb{F}_q^n$ via the linear transformation $\psi({x}) = {x}M$.

Next, we will show that \(\psi\) is a graph homomorphism. Let ${u}, {v}$ be adjacent vertices in $H(m, q, d)$ and ${w} = {u} - {v}$ be their difference vector. By the adjacency condition in the Hamming graph, ${w}$ contains exactly $d$ non-zero entries. Let ${z} = \psi({u}) - \psi({v}) = {w}M$ denote the difference vector in the target space $\mathbb{F}_q^n$. We verify that ${z}$ is a balanced vector.

Consider the first $(q-1)m$ coordinates of ${z}$, which are partitioned into $q-1$ blocks of length $m$. For each $k \in \{1, \dots, q-1\}$, the $k$-th block is given by $\lambda_k {w}$. If an entry $w_i$ is non-zero, it contributes the elements $\{ \lambda_1 w_i, \lambda_2 w_i, \dots, \lambda_{q-1} w_i \}$ to the vector ${z}$ across these blocks. Since $w_i \in \mathbb{F}_q^*$, the set $\{ \lambda w_i : \lambda \in \mathbb{F}_q^* \}$ is simply a permutation of $\mathbb{F}_q^*$. 
Because there are exactly $d$ such non-zero entries in ${w}$, each element of $\mathbb{F}_q^*$ appears exactly $d$ times in the first $(q-1)m$ coordinates of ${z}$. The total number of non-zero entries in ${z}$ is therefore $d(q-1)$. Since the total length of the vector is $n = qd$, the number of zero entries in ${z}$ is:
\[
    n - d(q-1) = qd - qd + d = d.
\]
Consequently, every element of $\mathbb{F}_q$ appears exactly $d$ times in ${z}$, which implies that ${z}$ is a balanced vector. 

By the definition of the graph $\varOmega^{(\mathbb{F}_q)}_{n}$, the vertices $\psi({u})$ and $\psi({v})$ are adjacent. Thus, $\psi$ is a graph homomorphism.
\end{proof}

\begin{proof}[Proof of \cref{thm:separate}]
    We first establish the lower bound for the classical chromatic number of $\varOmega^{(\mathbb{Z}_q)}_n$. Since $\chi(G)\alpha(G)\geq |V(G)|$, combined with \cref{lem:Zq ind-number}, we have 
    \[\chi(\varOmega^{(\mathbb{Z}_q)}_n)\geq \frac{q^n}{\alpha(G)}\geq(1+\varepsilon)^n\] for some $\varepsilon>0$.
    
    Next, we establish the lower bound for the classical chromatic number of $\varOmega^{(\mathbb{F}_q)}_n$. Frankl and R\"odl \cite[Theorem 1.10]{frankl1987forbidden} show that, for sufficiently large $n$, $\chi(H(n,q,d))\geq (1+\varepsilon)^n$.
    By \cref{lem: homo hamming to Fq}, we obtain 
    \[
    \chi(\varOmega^{(\mathbb{F}_q)}_{n})\geq \chi(H(\lfloor n/(q-1)\rfloor,q,d))\geq (1+\varepsilon)^{\lfloor n/(q-1)\rfloor}\geq (1+\varepsilon')^n,
    \]
    where $\varepsilon'>0$ is a suitable constant depending on $q$. This completes the proof.
\end{proof}

\section{Concluding remarks}\label{sec:rmk}
We investigate the quantum chromatic number of Hamming and generalized Hadamard graphs. Several interesting open questions arise from our work.

\paragraph{Orthogonal representations for $H(n,q,d)$ with $d < \frac{(q-1)n}{q}$} We develop a linear programming approach to construct modules-one orthogonal representations for Hamming graphs and establish upper bounds on $\xi'(H(n,q,d))$ for $d \le \frac{(q-1)n}{q}$. However, the bound in \cref{cor:OR-Case3} remains exponentially large.
\begin{question}
For $d = \delta n$ with $0 < \delta < \frac{q-1}{q}$, and with $d$ even when $q = 2$, is it always true that
\[
\xi'(H(n,q,d)) \le \operatorname{poly}(n),
\]
or does there exist some $d$ such that
\[
\xi'(H(n,q,d)) \ge \operatorname{exp}(n)?
\]
\end{question}

\paragraph{Exact value of $\chi_Q(H(n,q,\frac{(q-1)n}{q}))$ for $q \ge 3$.} We show that
\[
(q-1)n-(q-2) \le \chi_Q\big(H(n,q,\tfrac{(q-1)n}{q})\big) \le (q-1)n,
\]
but a gap of $(q-2)$ remains between the upper and lower bounds.
\begin{question}
How can this gap be closed?
\end{question}

\paragraph{Determining the minimum eigenvalue of $\varOmega_n^{(\mathbb{Z}_q)}$ for all $n$.} We determine the minimum eigenvalue of $\varOmega_n^{(\mathbb{Z}_q)}$ for sufficiently large $n$, and we conjecture that it remains the same for all feasible $n$. Formally, we state the conjecture as follows:
\begin{conjecture}
Let $q \ge 2$ be a positive integer, and let $n$ be divisible by $q$ such that $\frac{(q-1)n}{q}$ is an even integer. Then the minimum eigenvalue of $\varOmega_n^{(\mathbb{Z}_q)}$ is
\[
K^{(\mathbb{Z}_q)}_{\boldsymbol{n/q}}(n-2,1,0,\dots,0,1) = -\frac{\binom{n}{n/q,\ldots,n/q}}{n-1}.
\]
\end{conjecture}

\section*{Acknowledgements}
This research was conducted during the 2025 Quantum Information Summer School at Xidian University. The authors sincerely thank Professors Chong Shangguan and Sihuang Hu for their financial support as well as for their many valuable discussions and suggestions.

\section*{Declaration of generative AI and AI-assisted technologies in the manuscript preparation process}

During the preparation of this work, the authors used Gemini to improve the language and readability of the manuscript. After using this tool, the authors reviewed and edited the content as needed and take full responsibility for the content of the published article.

\bibliographystyle{plain}
\bibliography{references}
\end{document}